\documentclass[reqno]{amsart}

\usepackage{amsfonts}
\newtheorem{thm}{Theorem}[section]
\newtheorem{prop}{Proposition}[section]

\begin{document}
\newcounter{ZERO}
\catcode`\@=11

\@addtoreset{equation}{section}
\def\theequation{\ifnum\value{section}>0\relax
\thesection.\arabic{equation}\relax
\else\arabic{equation}\fi}

\newcommand{\C}{{{\mathbb C}}}
\newcommand{\R}{{{\mathbb R}}}
\newcommand{\Z}{{{\mathbb Z}}}
\newcommand{\bbC}{\mathbb{C}}
\newcommand{\bbR}{\mathbb{R}}
\newcommand{\bbH}{\mathbb{H}}
\newcommand{\bbN}{\mathbb{N}}
\newcommand{\bbZ}{\mathbb{Z}}
\newcommand{\mcP}{\mathcal{P}}
\newcommand{\mcH}{\mathcal{H}}
\newcommand{\mcS}{\mathcal{S}}
\newcommand{\mcE}{\mathcal{E}}
\newcommand{\mcR}{\mathcal{R}}
\newcommand{\la}{\lambda}
\newcommand{\al}{\alpha}
\newcommand{\be}{\beta}
\newcommand{\ga}{\gamma}
\newcommand{\de}{\delta}
\newcommand{\La}{\Lambda}
\newcommand{\gl}{\mathfrak{gl}}
\newcommand{\so}{\mathfrak{so}}
\renewcommand{\sp}{\mathfrak{sp}}
\newcommand{\F}[2]{F^{#1}_{(#2)}}
\newcommand{\E}[2]{E^{#1}_{(#2)}}
\newcommand{\V}[2]{V^{#1}_{(#2)}}
\newcommand{\wF}[2]{\widetilde F_{(#1)}^{#2}}
\newcommand{\wE}[2]{\widetilde E_{(#1)}^{#2}}
\newcommand{\wV}[2]{\widetilde V_{(#1)}^{#2}}
\renewcommand{\c}[3]{{c^{#1}_{{#2}\,{#3}}}}
\protect\newcommand{\ot}{\otimes}
\newcommand{\dual}[1]{ {\left({#1}\right)}^*}
\newcommand{\floor}[1]{\lfloor #1 \rfloor}
\newcommand{\nt}{{\noindent}}

\title[Stable branching rules]{Stable branching rules for\\ classical symmetric pairs}
\date{}
\author{Roger Howe, Eng-Chye Tan, Jeb F. Willenbring}

\maketitle

\section{ Introduction}

Given completely reducible representations, $V$ and $W$ of complex algebraic groups
$G$ and $H$ respectively, together with an embedding $H
\hookrightarrow G$, we let $[V,W] = \dim \text{Hom}_H \left(W,V
\right)$ where $V$ is regarded as a representation of $H$ by
restriction.  If $W$ is irreducible, then $[V,W]$ is the
multiplicity of $W$ in $V$.  This number may of course be infinite
if $V$ or $W$ is infinite dimensional.  A description of the
numbers $[V,W]$ is referred in the mathematics and physics
literature as a \emph{branching rule}.

The context of this paper has its origins in the work of D.
Littlewood.  In [Li2], Littlewood describes two classical
branching rules from a combinatorial perspective (see also [Li1]).
Specifically, Littlewood's results are branching multiplicities
for $GL_n$ to $O_n$ and $GL_{2n}$ to $Sp_{2n}$.  These pairs of
groups are significant in that they are examples of symmetric
pairs.  A symmetric pair is a pair of groups $(H,G)$ such that $G$
is a reductive algebraic group and $H$ is the fixed point set of a
regular involution defined on $G$.  It follows that $H$ is a
closed, reductive algebraic subgroup of $G$.

The goal of this paper is to put the formula into the context of
the first named author's theory of dual reductive pairs.  The
advantage of this point of view is that it relates branching from
one symmetric pair to another and as a consequence Littlewood's
formula may be generalized to all classical symmetric pairs.

Littlewood's result provides an expression for the branching
multiplicities in terms of the classical Littlewood-Richardson
coefficients (to be defined later) when the highest weight of the
representation of the general linear group lies in a certain
stable range.

The point of this paper is to show how when the problem of
determining branching multiplicities is put in the context of dual
pairs, a Littlewood-like formula results for any classical
symmetric pair.  To be precise, we consider 10 families of
symmetric pairs which we group into subsets determined by the
embedding of $H$ in $G$ (see Table I in \S 3).

\vskip 10pt

\nt {\bf Acknowledgement:} We are grateful to the referee for
elaborate and in-depth comments. Most notably, for the vast number
of references and comments on related works for branching rules and
tensor products of infinite-dimensional representations.

\subsection{Parametrization of Representations.} Let $G$ be a classical reductive algebraic group
over $\C$: $G = GL_n({\C})=GL_n$, the general linear group; or $G
= O_n ({\C})=O_n$, the orthogonal group; or $G =
Sp_{2n}({\C})=Sp_{2n}$, the symplectic group. We shall explain our
notations on irreducible representations of $G$ using integer
partitions.  In each of these cases, we select a Borel subalgebra
of the classical Lie algebra as is done in [GW].  Consequently,
all highest weights are parameterized in the standard way (see
[GW]).

A non-negative integer partition $\lambda$, with $k$ parts, is an
integer sequence $\lambda_1 \geq \lambda_2 \geq \ldots \geq
\lambda_k > 0$.  Sometimes we may refer to partitions as Young or
Ferrers diagrams.  We use the same notation for partitions as is
done in [Ma]. For example, we write $\ell(\la)$ to denote the
length (or depth) of a partition, $|\la|$ for the size of a
partition (i.e., $|\la| = \sum_i \la_i$). Also, $\la^\prime$
denotes the transpose (or conjugate) of $\la$ (i.e.,
$(\la^\prime)_i = |\{ \la_j: \la_j \geq i\}|$).  A partition where
all parts are even is called an {\it even} partition, and we shall
denote an even partition $2\delta_1 \geq 2\delta_2 \geq \ldots
\geq 2\delta_k$ simply by $2\delta$.

\bigskip
\noindent{\bf ${\bf GL_n}$ Representations:} Given non-negative integers
$p$ and $q$ such that $n \geq p+q$ and non-negative integer
partitions $\la^+$ and $\la^-$ with $p$ and $q$ parts
respectively, let $\F{(\la^+, \la^-)}{n}$ denote the irreducible
rational representation of $GL_n$ with highest weight given by the
$n$-tuple:
\[
    \begin{array}{ccc}
        (\la^+,\la^-) &=& \underbrace{\left(\la^+_1, \la^+_2, \cdots, \la^+_p, 0, \cdots, 0, -\la^-_{q},
\cdots, -\la^-_{1} \right)} \\
        & & n
    \end{array}
\]
If $\la^- = (0)$ then we will write $\F{\la^+}{n}$ for $\F{(\la^+,
\la^- )}{n}$. Note that if $\la^+ = (0)$ then
$\dual{\F{\la^-}{n}}$ is equivalent to $\F{(\la^+, \la^-)}{n}$.

\bigskip
\noindent{\bf ${\bf O_n}$ Representations:} The complex (or real) orthogonal group has two connected
components. Because the group is disconnected we cannot index
irreducible representation by highest weights. There is however an
analog of Schur-Weyl duality for the case of $O_n$ in which each
irreducible rational representation is indexed uniquely by a
non-negative integer partition $\nu$ such that $(\nu^\prime)_1 +
(\nu^\prime)_2 \leq n$.  That is, the sum of the first two columns
of the Young diagram of $\nu$ is at most $n$. (See [GW] Chapter 10
for details.)  Let $\E{\nu}{n}$ denote the irreducible
representation of $O_n$ indexed by $\nu$ in this way.

The irreducible rational representations of $SO_n$ may be indexed
by their highest weight, since the group is a connected reductive
linear algebraic group. In [GW] Section 5.2.2, the irreducible
representations of $O_n$ are determined in terms of their
restrictions to $SO_n$ (which is a normal subgroup having index
2). See [GW] Section 10.2.4 and 10.2.5 for the correspondence
between this parametrization and the above parametrization by
partitions.

\bigskip
\noindent{\bf ${\bf Sp_{2n}}$ Representations:} For a non-negative integer partition $\nu$ with $p$ parts where $p \leq
n$, let $\V{\nu}{2n}$ denote the irreducible rational
representation of $Sp_{2n}$ where the highest weight indexed by
the partition $\nu$ is given by the $n$ tuple:
\[
    \begin{array}{c}
        \underbrace{\left(\nu_1, \nu_2, \cdots, \nu_p, 0, \cdots, 0 \right)} \\ n
    \end{array}.
\]

\subsection{Littlewood-Richardson Coefficients}

Fix a positive integer $n_0$.  Let $\la$, $\mu$ and $\nu$ denote
non-negative integer partitions with at most $n_0$ parts.  For any
$n \geq n_0$ we have:
\[
\left[ \F{\mu}{n} \ot \F{\nu}{n}, \F{\la}{n} \right] = \left[ \F{\mu}{n_0}
\ot\F{\nu}{n_0}, \F{\la}{n_0} \right]
\]
And so we define:
\[
    \c{\la}{\mu}{\nu} := \left[ \F{\mu}{n} \ot \F{\nu}{n}, \F{\la}{n} \right]
\] for some (indeed any) $n \geq n_0$.

The numbers $\c{\la}{\mu}{\nu}$ are known as the
Littlewood-Richardson coefficients and are extensively studied in
the algebraic combinatorics literature.  Many treatments defined
from wildly different points of view.  See [BKW], [CGR], [Fu],
[GW], [JK], [Kn1], [Ma], [Sa], [St3] and [Su] for examples.

\subsection{Stability and the Littlewood Restriction Rules}
We now state the Littlewood restriction rules.

\begin{thm}[$O_n \subseteq GL_n$]
Given $\la$ such that $\ell(\la) \leq \frac{n}{2}$ and $\mu$ such
that $(\mu^\prime)_1 + (\mu^\prime)_2 \leq n$ then,
\begin{equation}
    [\F{\la}{n},\E{\mu}{n}] = \sum_{2\delta} c^\la_{\mu \; 2\delta}
\end{equation}
where the sum is over all non-negative even integer partitions
$2\delta$.
\end{thm}

\begin{thm}[$Sp_{2n} \subseteq GL_{2n}$]
Given $\la$ such that $\ell(\la) \leq n$ and $\mu$ such that
$\ell(\mu) \leq n$ then,
\begin{equation}
   [\F{\la}{2n}, \V{\mu}{2n}] = \sum_{2\delta} c^\la_{\mu \; (2\delta)^\prime}
\end{equation}
where the sum is over all non-negative integer partitions with even columns
$(2\delta)^\prime$.
\end{thm}

Notice that the hypothesis of the above two theorems do not
include an arbitrary parameter for the representation of the
general linear group.  The parameters which fall within this range
are said to be in the \emph{stable range}.  These hypothesis are
necessary but for certain $\mu$ it is possible to weaken them
considerably see [EW1] and [EW2].

    One purpose of this paper is to make the first steps toward a
uniform stable range valid for all symmetric pairs.  In the
situation presented here we approach the stable range on a
case-by-case basis.  Within the stable range, one can express the
branching multiplicity in terms of the Littlewood-Richardson
coefficients. These kinds of branching rules will later be
combined with the rich combinatorics literature on the
Littlewood-Richardson coefficients to provide more algebraic
structure to branching rules.

\section{Statement of the results.}

We now state our main theorem.  As it addresses 10 families of
symmetric pairs, which state in 10 parts.  The parts are grouped
into 4 subsets named:  Diagonal, Direct Sum, Polarization and
Bilinear Form.  These names describe the embedding of $H$ into $G$
(see Table I of \S 3).

\newpage

\nt {\bf Main Theorem:}

\subsection{Diagonal:}

\subsubsection{${\bf GL_n \subset GL_n \times GL_n}$}
Given non-negative integers, $p$, $q$, $r$ and $s$ with $n \geq
p+q+r+s$.  Let $\la^+$, $\mu^+$, $\nu^+$, $\la^-$, $\mu^-$,
$\nu^-$ be non-negative integer partitions.  If $\ell(\la^+) \leq
p+r$, $\ell(\la^-) \leq q+s$, $\ell(\mu^+) \leq p$, $\ell(\mu^-)
\leq q$, $\ell(\nu^+) \leq r$ and $\ell(\nu^-) \leq s$, then
\[
\left[ \F{(\mu^+, \mu^-)}{n} \ot \F{(\nu^+,\nu^-)}{n},
\F{(\la^+,\la^-)}{n} \right] = \sum \c{\la^+}{\al_2}{\al_1}
\c{\mu^+}{\al_1}{\ga_1} \c{\nu^-}{\ga_1}{\be_2}
\c{\la^-}{\be_2}{\be_1} \c{\mu^-}{\be_1}{\ga_2}
\c{\nu^+}{\ga_2}{\al_2}
\]
where the sum is over non-negative integer partitions $\al_1$,
$\al_2$, $\be_1$, $\be_2$, $\ga_1$ and $\ga_2$.

\subsubsection{${\bf O_n \subset O_n \times O_n}$} Given non-negative
integer partitions $\la$, $\mu$ and $\nu$ such that $\ell(\la )
\leq \floor{n/2}$ and $\ell(\mu) + \ell(\nu) \leq \floor{n/2}$, then
\[
\left[ \E{\mu}{n} \ot \E{\nu}{n}, \E{\la}{n} \right] = \sum \c{\la}{\al}{\be} \c{\mu}{\al}{\ga} \c{\nu}{\be}{\ga}
\]
where the sum is over all non-negative integer partitions $\al,
\be, \ga$.

\subsubsection{${\bf Sp_{2n} \subset Sp_{2n} \times Sp_{2n}}$} Given
non-negative integer partitions $\la$, $\mu$ and $\nu$ such that
$\ell(\la) \leq n$ and $\ell(\mu) + \ell(\nu) \leq n$, then
\[
\left[ \V{\mu}{2n} \ot \V{\nu}{2n}, \V{\la}{2n} \right] = \sum
\c{\la}{\al}{\be} \c{\mu}{\al}{\ga} \c{\nu}{\be}{\ga}
\]
where the sum is over all non-negative integer partitions $\al,
\be, \ga$.

\subsection{Direct Sum:}

\subsubsection{${\bf GL_n \times GL_m \subset GL_{n+m}}$}
Let $p$ and $q$ be non-negative integers such that $p+q \leq
\min(n,m)$.  Let $\la^+, \mu^+, \nu^+$ and $\la^-, \mu^-, \nu^-$
be non-negative integer partitions.  If $\ell(\la^+)$, $\ell(\mu^+)$, $\ell(\nu^+) \leq p $ and
$\ell(\la^-), \ell(\mu^-), \ell(\nu^-) \leq q $, then
\[
\left[ \F{(\la^+,\la^-)}{n+m}, \F{(\mu^+,\mu^-)}{n} \ot
\F{(\nu^+,\nu^-)}{m} \right] = \sum \c{\ga^+}{\mu^+}{\nu^+}
\c{\ga^-}{\mu^-}{\nu^-} \c{\la^+}{\ga^+}{\de}
\c{\la^-}{\ga^-}{\de}
\] where the sum is over all non-negative integer partitions $\ga^+$,
$\ga^-$, $\de$.

\subsubsection{${\bf O_n \times O_m \subset O_{n+m}}$}

Let $\la$, $\mu$ and $\nu$ be non-negative integer partitions such
that $\ell(\la), \ell(\mu), \ell(\nu) \leq \frac{1}{2} \min(n,m)$, then
\[
    \left[ \E{\la}{n+m}, \E{\mu}{n} \ot \E{\nu}{m} \right] =
    \sum \c{\ga}{\mu}{\nu} \c{\la}{\ga}{2\de}
\] where the sum is over all non-negative integer partitions $\de$ and
$\ga$.

\subsubsection{${\bf Sp_{2n} \times Sp_{2m} \subset Sp_{2(n+m)}}$}

Let $\la$, $\mu$ and $\nu$ be non-negative integer partitions such
that $\ell(\la), \ell(\mu), \ell(\nu) \leq \min(n,m)$, then
\[
    \left[ \V{\la}{2(n+m)}, \V{\mu}{2n} \ot \V{\nu}{2m} \right] =
    \sum \c{\ga}{\mu}{\nu} \c{\la}{\ga}{(2\de)^\prime}
\] where the sum is over all non-negative integer partitions $\de$ and
$\ga$.

\subsection{Polarization:}

\subsubsection{${\bf GL_n \subset O_{2n}}$}
Let $\mu^+$, $\mu^-$ and $\la$ be non-negative integer partitions
with at most $\floor{n/2}$ parts, then
\[
\left[ \E{\la}{2n},\F{(\mu^+,\mu^-)}{n} \right] = \sum \c{\ga}{\mu^+}{\mu^-}
\c{\la}{\ga}{(2\de)^\prime}
\]
where the sum is over all non-negative integer partitions $\de$
and $\ga$.

\subsubsection{${\bf GL_n \subset Sp_{2n}}$}
Let $\mu^+$, $\mu^-$ and $\la$ be non-negative integer partitions
with at most $\floor{n/2}$ parts, then
\[
\left[ \V{\la}{2n},\F{(\mu^+,\mu^-)}{n} \right] = \sum \c{\ga}{\mu^+}{\mu^-}
\c{\la}{\ga}{2\de}
\]
where the sum is over all non-negative integer partitions $\de$
and $\ga$.

\subsection{Bilinear Form:}

\subsubsection{${\bf O_n \subset GL_n}$}
Let $\la^+$, $\la^-$ and $\mu$ denote non-negative integer
partitions with at most $\floor{n/2}$ parts, then
\[
    \left[ \F{(\la^+,\la^-)}{n}, \E{\mu}{n} \right]
    = \sum \c{\mu}{\al}{\be} \c{\la^+}{\al}{2\ga} \c{\la^-}{\beta}{2\de}
\] where the sum is over all non-negative integer partitions
$\al$, $\beta$, $\ga$ and $\de$.

\subsubsection{${\bf Sp_{2n} \subset GL_{2n}}$}
Let $\la^+$, $\la^-$ and $\mu$ denote non-negative integer
partitions with at most $n$ parts, then
\[
    \left[ \F{(\la^+,\la^-)}{2n}, \V{\mu}{2n} \right]
    = \sum \c{\mu}{\al}{\be} \c{\la^+}{\al}{(2\ga)^\prime} \c{\la^-}{\beta}{(2\de)^\prime}
\] where the sum is over all non-negative integer partitions
$\al$, $\beta$, $\ga$ and $\de$.

\bigskip
\subsection{\bf Remarks:} Although a thorough survey is beyond our
present goals, we wish to record here many previous works on
branching rules which in many cases overlaps with ours.  We are
grateful to the referee who has given us an extensive list of
references with comments on related works by experts. We shall
briefly summarize related works as follows:

\begin{enumerate}

\item[(a)] {\bf Diagonal:} The first rule 2.1.1 appears as (4.6)
with (4.15) in King's paper [Ki2]. The branching rules 2.1.2 and
2.1.3 for orthogonal and symplectic groups goes back to Newell
[Ne] and Littlewood [Li3].  A more rigorous account of the
Diagonal rules also appears in [Ki4], along with a treatment of
rational representations of $GL_n$.  See Theorem 4.5 and Theorem
4.1 of [Ki4] and the references therein.  Our methods are also
cast in this same generality. Further, 2.1.2 and 2.1.3 are
beautifully presented in Sundaram's survey [Su] with references to
the proofs in [BKW].

\item[(b)] {\bf Direct Sum:} Rule 2.2.1 appears as (5.8) with
(4.16) in one of the earlier works of King [Ki1], which derives
from a conjecture in the Ph.D. thesis by Abramsky [Ab].  These
branching rules are also addressed in [Ko] and [KT].
Specifically, 2.2.1 can be found in Proposition 2.6 of [Ko], and
2.2.2 and 2.2.3 can be found in Theorem 2.5 and Corollary 2.6 of
[KT].  An account of the Direct Sum rules also appears in [Ki4]
(see (2.1.6) and the references therein).

\item[(c)] {\bf Polarization:} The polarization branching rules,
2.3.1 and 2.3.2 are stated as (4.21) and (4.22), respectively, in
[Ki3], and also as Theorem A1 of [KT].

\item[(d)] {\bf Bilinear Form:} The Littlewood restriction rule is
a special case of formulas, 2.4.1 and 2.4.2 (see [Li1] and [Li2]).
These two formulas can be viewed as a generalization of
Littlewood's restriction rule. Besides the Diagonal branching
formulas, [Su] also presents a thorough treatment of the classical
Littlewood restriction rules.  However, in the most general form,
the rules 2.4.1 and 2.4.2 appear as (5.7) with (4.19), and (5.8)
with (4.23) respectively in [Ki2].

\end{enumerate}

Most of the results have been sufficiently well known by experts.
For a well presented survey of the representation theory of the
classical groups from a combinatorial point of view we refer the
reader to [Su].  And the late Wybourne and his students have even
incorporated these results in the software package SCHUR
downloadable at http://smc.vnet.net/Schur.html.  This package
implements all the modification rules given in [Ki2] and [BKW]
that allow the stable branching rules to be generalized so as to
cover all possible non-stable cases as well.

From our point of view, it is striking that the theory of dual
pairs leads to proofs of all 10 of these formulas in such a
unified manner.  We feel that this unifying theme should be
brought out in the literature more systematically than it has
been.

In [Su] Theorem 5.4, it is shown how the Littlewood Richardson
rule for branching from $GL_{2n}$ to $Sp_{2n}$ may be modified to
obtain a version of the Littlewood restriction rule which is valid
outside the stable range. Removing the stability condition for
Littlewood's restriction rules is a delicate problem, which was
also addressed in [EW1] and [EW2]. Classically, Newell [Ne]
presents modification rules to the Littlewood restriction rules to
solve the branching problem outside of the stable range (see [Su]
and [Ki2]). For some recent remarks on the literature of branching
rules we refer the reader to [Ki3], [Kn2], [Kn3], [Kn4] and [Pr].
The discussions in [Kn2] are relevant to our approach.  Some of
the results in [Kn2] are important special cases of the results in
[GK].

While we only require decompositions of tensor products of
infinite- \linebreak dimensional holomorphic discrete series as in
[Re], we also wish to note the numerous works in more generality
which can be found in [RWB], [KW], [TTW], [OZ] and [KTW], amongst
the many others. It is interesting to note that the papers [RWB],
[KW], [TTW] and [KTW] have exploited the duality correspondence
(see \S 3.1) to relate the multiplicities of tensor products of
infinite-dimensional representations to multiplicities of tensor
products of finite-dimensional representations in the same spirit
of [Ho1].

\section{Dual Pairs and Reciprocity}

The formulation of
classical invariant theory in terms of dual pairs \cite{howe-remarks} allows one to
realize branching properties for classical symmetric pairs by considering concrete
realizations of representations on algebras of polynomials on vector spaces.

\subsection{Dual Pairs and Duality Correspondence} Let $W \simeq \R^{2m}$ be a\linebreak $2m$-dimensional real vector space with symplectic form $<\cdot,\cdot>$. Let $Sp(W)=Sp_{2m}(\R)$ denote the isometry group of the form $<\cdot,\cdot>$.  A pair of subgroups $(G,G')$ of $Sp_{2m}(\R)$ is called a {\it reductive dual pair} (in $Sp_{2m}(\R)$) if
\begin{enumerate}
\item[(a)] $G$ is the centralizer of $G'$ in $Sp_{2m}(\R)$ and vice versa, and
\item[(b)] both $G$ and $G'$ act reductively on $W$.
\end{enumerate}

The fundamental group of $Sp_{2m}(\R)$ is the fundamental group of
$U_m$, its maximal compact subgroup, and is isomorphic to $\Z$.
Let $\widetilde{Sp}_{2m}(\R)$ denote a choice of a double cover of
$Sp_{2m}(\R)$.  We will refer to this as the {\it metaplectic
group}. Also let $\widetilde U_m$ denote the pull-back of the
covering map on $U_m$.  Shale-Weil constructed a distinguished
representation $\omega$ of $\widetilde{Sp}_{2m}(\R)$, which we
shall refer to as the {\it oscillator representation}.  This is a
unitary representation and one realization is on the space of
holomorphic functions on $\C^m$, commonly referred to as the Fock
space. In this realization, the $\widetilde{U}_m$-finite functions
appear as polynomials on $\C^m$ which we denote as $\mcP(\C^m)$. A
vector $v \in \mcP(\C^m)$ is $\widetilde{U}_m$-finite if the span
of $\widetilde{U}_m \cdot v$ in $\mcP(\C^m)$ is finite
dimensional.

Choose $z_1,z_2,\ldots, z_m$ as a system of coordinates on $\C^m$.  The Lie algebra action of $\sp_{2m}$ (the complexified Lie algebra of $Sp_{2m}(\R)$) on $\mcP(\C^m)$ can be described by the following operators:
\begin{equation}
\omega (\sp_{2m}) = \sp_{2m}^{(1,1)} \oplus \sp_{2m}^{(2,0)}\oplus  \sp_{2m}^{(0,2)}
\end{equation}
where
\begin{equation}
\begin{aligned}
\sp_{2m}^{(1,1)}  &= \text{ Span }\left\{ \frac{1}{2} \left( z_i \frac{\partial}
{\partial z_j} + \frac{\partial}{\partial z_i} z_j \right) \right\},\\
\qquad\ \qquad \qquad \qquad \sp_{2m}^{(2,0)}  &= \text{ Span }\{ z_iz_j  \},\qquad \qquad \qquad \qquad \qquad \qquad \qquad \ \\
\sp_{2m}^{(0,2)}  &= \text{ Span }\left\{ \frac{\partial^2}{\partial z_i \partial z_j} \right\}.
\end{aligned}
\end{equation}

The decomposition (3.1) is an instance of the complexified Cartan decomposition
\begin{equation}
\sp_{2m} = \frak k \oplus \frak {p}^{+} \oplus \frak {p}^{-}
\end{equation}
where $\sp_{2m}^{(1,1)}\simeq \omega (\frak k)$, $\sp_{2m}^{(2,0)} \simeq \omega (\frak {p}^+)$ and $\sp_{2m}^{(0,2)} \simeq \omega (\frak {p}^-)$.
If $\mcP(\C^m) = \sum_{s\geq 0} \mcP^s (\C^m)$ is the natural grading on $\mcP(\C^m)$, it is immediate that $\sp_{2m}^{(i,j)}$ brings $\mcP^s (\C^m)$ to $\mcP^{s+i-j} (\C^m)$.

Let us restrict our dual pairs to the following:
\begin{equation}
(O_n(\R), Sp_{2k}(\R)),\qquad (U_n,U_{p,q}),\qquad (Sp(n),O^\ast_{2k}).
\end{equation}
Observe that the first member is compact, and these pairs are usually loosely referred to as {\it compact pairs}.

To avoid technicalities involving covering groups, instead of the
real groups $(G_0, G_0')$, we shall discuss in the context of
pairs $(G,\frak g')$ where $G$ is a complexification of $G_0$ and
$\frak g'$ is a complexification of the Lie algebra of $G_0'$. The
use of the phrase ``up to a central character'' in the statements
(a) to (c) below basically suppresses the technicalities involving
covering groups. Each of these pairs can be conveniently realized
as follows:
\begin{enumerate}
\item[(a)] $(O_n(\R), Sp_{2k}(\R)) \subset Sp_{2nk}(\R)$:\newline
Let $\C^n \otimes \C^k$ be the space of $n$ by $k$ complex matrices. The complexified pair $(O_n, \sp_{2k})$ acts on $\mcP(\C^n \otimes \C^k)$ which are the $\widetilde{U}_{nk}$-finite functions. The group $O_n$ acts by left multiplication on $\mcP(\C^n \otimes \C^k)$ and can be identified with the holomorphic extension of the $O_n(\R)$ action on the Fock space.  The action of the subalgebra $\gl_k$ of $\sp_{2k}$ is (up to a central character) the derived action coming from the natural right action of multiplication by $GL_k$.
\vskip 10pt

\item[(b)] $(U_n, U_{p,q}) \subset Sp_{2n(p+q)}(\R)$:\newline
For this pair, we may identify the $\widetilde{U}_{n(p+q)}$-finite functions with the polynomial ring $\mcP(\C^n \otimes \C^p \oplus (\C^{n})^{\ast} \otimes \C^q)$. The complexified pair is $(GL_n, \gl_{p,q})$.  There is a natural action of $GL_n$ and $GL_p \times GL_q$ on this polynomial ring as follows:
$$
(g,h_1,h_2)\cdot F(X,Y) = F(g^{-1}Xh_1,g^tYh_2),
$$
where $X\in \C^n \otimes \C^p$, $Y\in (\C^n)^\ast \otimes \C^q$, $g\in GL_n$, $h_1 \in GL_p$ and $h_2 \in GL_q$. Obviously both left and right actions commute.  Here $\gl_{p,q} \simeq \gl_{p+q}$, but we choose to differentiate the two because of the role of the subalgebra $\gl_p \oplus \gl_q$, which acts by (up to a central character) the derived action of $GL_p \times GL_q$ on the polynomial ring.
\vskip 10pt

\item[(c)] $(Sp(n), O^*_{2k}) \subset Sp_{4nk}(\R)$:\newline
In this case, $\mcP(\C^{2n} \otimes \C^k)$ are the $\widetilde{U}_{2nk}$-finite functions, with natural left and right actions by $Sp_{2n}$ and $GL_k$ respectively. The complexified pair is $(Sp_{2n}, \frak o_{2k})$, where the subalgebra $\gl_k$ of $\frak o_{2k}$ acts by (up to a central character) the derived right action of $GL_k$.

\end{enumerate}

With the realizations of these compact pairs $(G,\frak g') \subset Sp_{2m}(\R)$, let us look at the representations that appear. Form
$$
\frak {g'}^{(i,j)} = \sp_{2m}^{(i,j)} \cap \omega (\frak g')
$$
to get
\begin{equation}
\omega (\frak g') = \frak {g'}^{(1,1)} \oplus \frak {g'}^{(2,0)} \oplus \frak {g'}^{(0,2)}.
\end{equation}

Observe that $G_0'$ is Hermitian symmetric in all three cases, and the decomposition above is an instance of the complexified Cartan decomposition
\begin{equation}
\frak g' = \frak k' \oplus \frak {p'}^{+} \oplus \frak {p'}^{-}
\end{equation}
where $\frak {g'}^{(1,1)}\simeq \omega (\frak k')$, $\frak {g'}^{(2,0)} \simeq \omega (\frak {p'}^+)$ and $\frak {g'}^{(0,2)} \simeq \omega (\frak {p'}^-)$.  In particular, $\frak k'$ has a one-dimensional center and $\frak {p'}^{\pm}$ are the $\pm i$ eigenspaces of this center.  Each $\frak {p'}^{\pm}$ is an abelian Lie algebra.  Note, in particular,
\begin{equation}
[\frak {g'}^{(1,1)},\frak {g'}^{(2,0)}] \subset \frak {g'}^{(2,0)} \qquad \text{and} \qquad
 [\frak {g'}^{(1,1)},\frak {g'}^{(0,2)}] \subset \frak {g'}^{(0,2)}.
\end{equation}

A representation $(\rho, V_\rho)$ of $\frak g'$ is {\it holomorphic} if there is a non-zero vector $v_0 \in V_\rho$ killed by $\rho (\frak {p'}^-)$.  The following are the key properties of holomorphic representations:
\begin{enumerate}
\item[(a)] There is a non-trivial subspace
$$
(V_\rho )_0 =\ker \rho (\frak {p'}^{-}) = \{ v \in V_\rho \mid \rho (Y)\cdot v=0 \text{ for all } Y \in \frak{p'}^-   \}
$$
which is $\frak k'$ irreducible. This is known as the {\it lowest $\frak k'$-type} of $\rho$.

\item[(b)] $V_\rho$ is generated by $(V_\rho)_0$, more precisely,
$$
V_\rho= \mathcal U (\frak {p'}^+) \cdot V_0 = \mcS (\frak {p'}^+) \cdot V_0.
$$
The second equality results because $\frak {p'}^+$ is abelian.

\end{enumerate}

Now one of the key features in the formalism of dual pairs is the branching decomposition of the oscillator representation.  The branching property for compact pairs  alluded to is (see \cite{howe-remarks} and the references therein):
\begin{equation}
\mcP(\C^m) \mid_{G \times \frak g'} = \bigoplus_{\tau \in S \subset \widehat{G}} \tau \otimes V_{\tau'}
\end{equation}
where $S$ is a subset of the set of irreducible representations of $G$, denoted by $\widehat{G}$.  The representations $V_{\tau'}$ (written to emphasize the correspondence $\tau \leftrightarrow \tau'$ and the dependence on $\tau \in \widehat{G}$) are irreducible holomorphic representations of $\frak g'$.  They are known to be derived modules of irreducible unitary representations of some appropriate cover of $G_0'$ (\cite{howe-transcending}).  The key feature of this branching is the uniqueness of the correspondence, i.e., a representation of $G$ appearing uniquely determines the representation of the $\frak g'$ module that appears and vice-versa. We refer to this as the {\it duality correspondence}.

This duality is subjugated to another correspondence in the space of harmonics.

\begin{thm}(\cite{howe-remarks}, \cite{kashiwara-vergne})
Let $\mathcal H = \ker \frak {g'}^{(0,2)}$ be the space of harmonics. Then $\mathcal H$ is a $G \times K'$ module and it admits a multiplicity-free $G\times K'$ (hence $G \times \frak k'$) decomposition:
\begin{equation}
\mathcal H = \bigoplus_{\tau \in S \subset \widehat{G}} \tau \otimes \ker \rho_{\tau'}( \frak {g'}^{(0,2)}).
\end{equation}
We also have the separation of variables theorem providing the following $G \times \frak g'$ decomposition:
\begin{equation}
\begin{aligned}
\qquad \mcP(\C^m) &= \mathcal H \cdot \mcS (\frak {g'}^{(2,0)} )  = \left\{ \bigoplus_{\tau \in S \subset \widehat{G'}} \tau \otimes \ker \rho_{\tau'} ( \frak {g'}^{(0,2)}) \right\}\cdot \mcS (\frak {g'}^{(2,0)} ) \\
 &= \bigoplus_{\tau \in S \subset \widehat{G}} \tau \otimes \left\{ \mcS (\frak {g'}^{(2,0)} ) \cdot \ker \rho_{\tau'} ( \frak {g'}^{(0,2)})  \right\}  =  \bigoplus_{\tau \in S \subset \widehat{G}} \tau \otimes V_{\tau'} .
\end{aligned}
\end{equation}

\end{thm}

The structure of $V_{\tau^\prime}$ is even nicer in certain category of pairs, which we will refer to as the {\it stable range}. The stable range refers to the following:
\begin{enumerate}
\item[(a)] $(O_n(\R),Sp_{2k}(\R))$ for $n \geq 2k$;
\item[(b)] $(U_n,U_{p,q})$ for $n \geq p+q$;
\item[(c)] $(Sp(n), O^*_{2k})$ for $n \geq k$.
\end{enumerate}

\noindent In the stable range, the holomorphic representations of
$\frak g'$ that occur have $\frak k'$-structure which are nicer
(\cite{harish-chandra}, \cite{schmid1}, \cite{schmid2}), namely,
$$
V_{\tau'} = \mcS (\frak {g'}^{(2,0)} ) \otimes \ker \tau' ( \frak {g'}^{(0,2)}).
$$
They are known as {\it holomorphic discrete series} or {\it limits of holomorphic discrete series} (in some limiting cases of the parameters determining $\tau^\prime$) of the appropriate covering group of $G_0'$.  It is these representations that will feature prominently in this paper.

Let us conclude by describing the duality correspondence for the
compact dual pairs in the stable range.  Parts of the following
well known result can be found in several places, see \cite{ehw},
\cite{harish-chandra}, \cite{howe-remarks},
\cite{howe-transcending}, \cite{schmid1}, \cite{schmid2} for
example.

\begin{thm}
\begin{enumerate}
\item[(a)] $(O_n(\R) ,Sp_{2k}(\R))$: The duality correspondence for $O_n \times \sp_{2k}$ is:
\begin{equation}
\mcP(\C^n \otimes \C^k) = \bigoplus_{\lambda} E_{(n)}^\lambda \ot
\widetilde{E}^\lambda_{(2k)}
\end{equation}
where $\lambda$ runs through the set of all non-negative integer
partitions such that $l(\lambda) \leq k$ and $(\lambda')_1 +
(\lambda')_2 \leq n$. The space $\widetilde{E}^\lambda_{(2k)}$ is
an irreducible holomorphic representation of $\sp_{2k}$ of lowest
$\gl_k$-type $F_{(k)}^\lambda$.  In the stable range $n \geq 2k$,
$$
\widetilde{E}^\lambda_{(2k)} \simeq \mcS (\sp_{2k}^{(2,0)})
\otimes F_{(k)}^\la \simeq \mcS (\mcS^2 \C^k) \otimes F_{(k)}^\la
.
$$

\vskip 10pt

\item[(b)] $(U_n,U_{p,q})$: The duality correspondence for $GL_n \times \gl_{p,q}$ is:
\begin{equation}
\mcP(\C^n \otimes \C^p \otimes (\C^n)^\ast \otimes \C^q) =
\bigoplus_{\lambda^+, \lambda^-} F_{(n)}^{(\lambda^+,\lambda^-)}
\ot  \widetilde{F}^{(\lambda^+,\lambda^-)}_{(p,q)}
\end{equation}
where the sum is over all non-negative integer partitions
$\lambda^+$ and $\lambda^-$ such that $l(\lambda^+) \leq p$,
$l(\lambda^-) \leq q$ and $l(\lambda^+) + l(\lambda^-) \leq n$.
The space $\widetilde{F}^{(\lambda^+,\lambda^-)}_{(p,q)}$ is an
irreducible holomorphic representation of $\gl_{p,q}$ with lowest
$\gl_p \oplus \gl_q$-type $F_{(p)}^{\lambda^+} \ot
F_{(q)}^{\lambda^-}$. In the stable range $n \geq p+q$,
$$
\widetilde{F}^{(\lambda^+,\lambda^-)}_{(p,q)} \simeq \mcS
(\gl_{p,q}^{(2,0)}) \otimes F_{(p)}^{\lambda^+} \ot
F_{(q)}^{\lambda^-} \simeq \mcS (\C^p \otimes \C^q) \otimes
F_{(p)}^{\lambda^+} \ot F_{(q)}^{\lambda^-}.
$$

\nt The degenerate case when $q=0$ is particularly interesting.
This is the $GL_n \times GL_p$ duality:
\begin{equation}
\mcP(\C^n \otimes \C^p ) =
\bigoplus_{\la} F_{(n)}^{\la} \ot F_{(p)}^{\la}
\end{equation}
where the sum is over all integer partitions $\lambda$ such that $l(\lambda) \leq \min (n,p)$.

\vskip 10pt \item[(c)] $(Sp(n), O^*_{2k})$: The duality
correspondence for $Sp_{2n} \times {\so}_{2k}$ is:
\begin{equation}
\mcP(\C^{2n} \otimes \C^k) = \bigoplus_{\lambda} V_{(2n)}^\lambda
\ot \widetilde{V}^\lambda_{(2k)}
\end{equation}
where $\lambda$ runs through the set of all non-negative integer
partitions such that $l(\lambda) \leq \min (n,k)$.  The space
$\widetilde{V}^\lambda_{(2k)}$ is an irreducible holomorphic
representation of ${\so}_{2k}$ with lowest $\gl_k$-type
$F_{(k)}^\lambda$.  In the stable range $n \geq k$,
$$
\widetilde{V}^\lambda_{(2k)} \simeq \mcS ({\so}_{2k}^{(2,0)})
\otimes F_{(k)}^\la \simeq \mcS (\wedge^2 \C^k) \otimes
F_{(k)}^\la .
$$

\end{enumerate}
\end{thm}
\vskip 10pt

\subsection{Symmetric Pairs and Reciprocity Pairs} In the context of dual pairs, we would like to understand the branching of irreducible representations from $G$ to $H$, for symmetric pairs $(H, G)$.
Table I lists the symmetric pairs which we will cover in this paper.
\vskip 10pt

\begin{center} {\bf Table I: Classical Symmetric Pairs} \end{center}
\vskip 10pt
\begin{center}
\begin{tabular}{|c|c|c|} \hline
{\bf Description} & ${\bf H}$ & ${\bf G}$ \\ \hline
Diagonal & $GL_n$ & $GL_n \times GL_n$ \\ \hline
Diagonal & $O_n$ & $O_n \times O_n$ \\ \hline
Diagonal & $Sp_{2n}$ & $Sp_{2n} \times Sp_{2n} $ \\ \hline
Direct Sum & $GL_n \times GL_m$ & $GL_{n+m}$ \\ \hline
Direct Sum & $O_n \times O_m$ & $O_{n+m}$ \\ \hline
Direct Sum & $Sp_{2n} \times Sp_{2m}$ & $Sp_{2(n+m)}$ \\ \hline
Polarization & $GL_n$ & $O_{2n}$ \\ \hline
Polarization & $GL_n$ & $Sp_{2n}$ \\ \hline
Bilinear Form & $O_n$   & $GL_n$ \\ \hline
Bilinear Form & $Sp_{2n}$  & $GL_{2n}$ \\ \hline
\end{tabular}
\end{center}
\vskip 10pt

If $G$ is a classical  group over $\C$, then $G$ can be embedded as
one member of a dual pair in the symplectic group as described
in \cite{howe-remarks}. The resulting pairs of groups are
$(GL_n , GL_m )$ or $(O_n , Sp_{2m})$, each inside
$Sp_{2nm}$, and are called {\it irreducible} dual pairs. In general, a dual pair of reductive groups in $Sp_{2r}$ is a product of such pairs.

\vskip 10pt

\begin{prop} Let $G$ be a classical group, or a product of two copies
of a classical group. Let $G$ belong to a dual pair $(G, G')$ in a symplectic
group $Sp_{2m}$.  Let $H \subset G$ be a symmetric subgroup, and let $H'$ be
the centralizer of $H$ in $Sp_{2m}$. Then $(H, H')$ is also a dual pair in
$Sp_{2m}$, and $G'$ is a symmetric subgroup inside $H'$.
\end{prop}

\vskip 10pt

\nt {\bf Proof:} This can be shown by fairly easy case-by-case
checking. These are shown in Table II. We call these pair of pairs
{\it reciprocity pairs}.  These are special cases of see-saw pairs
[Ku].  \qquad $\square$

\vskip 10pt {\small
\begin{center} {\bf Table II: Reciprocity Pairs} \end{center}
\vskip 10pt
\begin{center}
\begin{tabular}{|c|c|c|} \hline
\ {\bf Symmetric Pair } ${\bf (H,G)}$ & ${\bf (H,\frak h')}$ &
${\bf (G,\frak g')}$ \\ \hline \hline \ $(GL_n, GL_n \times GL_n)$
\ &\ $(GL_n, \gl_{m+\ell})$ \ &\ $(GL_n \times GL_n, \gl_m \times
\gl_\ell)$ \ \\ \hline \ $(O_n, O_n \times O_n)$ \ &\ $(O_n,
\sp_{2(m+\ell )})$ \ &\ $(O_n \times O_n, \sp_{2m} \times
\sp_{2\ell})$ \ \\ \hline \ $(Sp_{2n}, Sp_{2n} \times Sp_{2n})$ \
&\ $(Sp_{2n}, \so_{2(m+\ell)})$ \ &\ $(Sp_{2n} \times Sp_{2n},
\so_{2m} \oplus  \so_{2\ell})$ \ \\ \hline \ $(GL_n \times GL_m,
GL_{n+m})$ \ &\ $(GL_n \times GL_m, \gl_\ell \times \gl_\ell)$ \
&\ $(GL_{n+m}, \gl_\ell )$ \ \\ \hline \ $(O_n \times O_m,
O_{n+m})$ \ &\ $(O_n \times O_m, \sp_{2\ell} \oplus \sp_{2\ell})$
\ &\ $(O_{n+m}, \sp_{2\ell})$ \ \\ \hline \ $(Sp_{2n} \times
Sp_{2m}, Sp_{2(n+m)})$ \ &\ $(Sp_{2n} \times Sp_{2m}, \so_{2\ell }
\oplus \so_{2\ell })$  \ &\ $(Sp_{2(n+m)}, \so_{2\ell })$ \ \\
\hline \ $(GL_n, O_{2n})$ \ &\ $(GL_n, \gl_{2m})$ \ &\ $(O_{2n},
\sp_{2m})$ \ \\ \hline \ $(GL_n, Sp_{2n})$ \ &\ $(GL_n, \gl_{2m})$
\ &\ $(Sp_{2n}, \so_{2m})$ \ \\  \hline \quad $(O_n, \gl_{n})$ \
&\ $(O_n, \sp_{2m})$ \ &\ $(GL_{n}, \gl_m)$ \ \\ \hline \
$(Sp_{2n}, GL_{2n})$ \ &\ $(Sp_{2n}, \so_{2m})$ \ &\ $(GL_{2n},
\gl_{m})$ \quad \\ \hline
\end{tabular}
\end{center}}

\vskip 10pt

Consider the dual pairs $(G,G')$ and $(H,H')$ in $Sp_{2m}$. We illustrate them in the see-saw manner as follows:
$$
\begin{aligned}
G &\qquad - \qquad G' \\
\cup &\qquad \quad\ \qquad \cap \\
H &\qquad - \qquad H'
\end{aligned}
$$
Recall the duality correspondence for $G \times \frak g'$ and $H \times \frak h'$ on the space $\mcP(\C^m)$:
$$
\begin{aligned}
\mcP(\C^m) \mid_{ G \times \frak g'} &= \bigoplus_{\sigma \in S \subset \widehat G} \sigma \otimes V_{\sigma'} =  \bigoplus_{\sigma \in S \subset \widehat G} \mcP(\C^m)_{\sigma \otimes \sigma'}, \\
\mcP(\C^m) \mid_{ H \times \frak h'} &= \bigoplus_{\tau \in T \subset \widehat H} \tau \otimes W_{\tau'} =  \bigoplus_{\tau \in T \subset \widehat H} \mcP(\C^m)_{\tau \otimes \tau'},
\end{aligned}
$$
where we have written $\mcP(\C^m)_{\sigma \otimes \sigma'}$ and
$\mcP(\C^m)_{\tau \otimes \tau'}$ as the $\sigma \otimes \sigma'$-isotypic component and
$\tau \otimes \tau'$-isotypic component in $\mcP(\C^m)$ respectively.  Given $\sigma$ and $\tau'$, we can seek the $\sigma \otimes \tau'$-isotypic component in $\mcP(\C^m)$ in two ways as follows:
\begin{equation}
\mcP(\C^m)_{\sigma \otimes \tau'}  \simeq (\sigma\mid_H )_\tau
\otimes V_{\sigma'} \simeq \tau \otimes ( W_{\tau'} \mid_{\frak
h'})_{\sigma'}. \end{equation} In other words, we have the
equality of multiplicities (as pointed out in [Ho1]):
\begin{equation}
[\sigma, \tau ] = [W_{\tau'} , V_{\sigma'}]
\end{equation}
that is, the multiplicity of $\tau$ in $\sigma\mid_H$ is equal to the multiplicity of $V_{\sigma'}$ in
$W_{\tau'} \mid_{\frak h'}$.  This is good enough for our purposes in this paper.  However, this equality of multiplicities is just a feature of some deeper phenomenon -- an isomorphism of certain branching algebras which captures the respective branching properties.

\subsection{Branching Algebras} One approach to branching problems exploits the fact that the
representations have a natural product structure, embodied by the algebra of regular functions on the flag manifold of the group. For a reductive complex algebraic $G$,
let $N_G$ be a maximal unipotent subgroup of $G$. The group $N_G$ is
determined up to conjugacy in $G$.  Let $A_G$ denote a maximal torus which normalizes $N_G$, so that $B_G = A_G\cdot N_G$ is a Borel subgroup of $G$. Let $\widehat{A}_G^+$ be the set
of dominant characters of $A_G$ -- the semigroup of highest weights of
irreducible representations of $G$.  It is well-known (see for instance, \cite{howe-schur}) that the space of regular functions on the
coset space  $G/N_G$ decomposes (under the action of $G$ by left translations) as a direct sum of one copy of each irreducible representation $V_{\psi}$, of highest weight $\psi$, of $G$:
$$
\mcR (G/N_G) \simeq \bigoplus_{ \psi \in \widehat{A}_G^+} V_{\psi}.
$$

We note that $\mcR (G/N_G)$ has the structure of an $\widehat{A}_G^+$-graded algebra, for which the $V_{\psi}$ are the graded components.  Let $H \subset G$ be a reductive subgroup and $A_H=A_G \cap H$ be a maximal torus of $H$ normalizing $N_H$, a maximal unipotent subgroup of $H$, so that $B_H=A_H \cdot N_H$ is a Borel subgroup of $H$. We consider the algebra $\mcR (G/N_G)^{N_H}$, of
functions on $G/N_G$ which are invariant under left translations by $N_H$ . This is an ($\widehat{A}_G^+ \times \widehat{A}_H^+$)-graded algebra. Knowledge of
$\mcR (G/N_G)^{N_H}$ as a ($\widehat{A}_G^+ \times \widehat{A}_H^+$)-graded algebra tell us
how representations of $G$ decompose when restricted to $H$, in other words,
it describes the branching rule from $G$ to $H$. We will call
$\mcR (G/N_G)^{N_H}$ the $(G, H)$ \it branching algebra. \rm  When $G \simeq H
\times H$, and $H$ is embedded diagonally in $G$, the branching algebra
describes the decomposition of tensor products of representations of $H$,
and we then call it the \it tensor product algebra \rm for $H$.

Let us explain briefly how branching algebras, dual pairs and reciprocity are related. For a reciprocity pair $(G,\frak g')$, $(H,\frak h')$, the duality correspondences are subjugated to a correspondence in the space of harmonics $\mcH$ (see Theorem 3.1).  Branching from holomorphic discrete series of $\frak h'$ to $\frak g'$ behaves very much like finite-dimensional representations in relation to their highest weights and is captured entirely by the branching from the lowest $K_{H'}$-type to $K_{G'}$. Although $\mcH$ is not an algebra, it can still be identified as a quotient algebra of $\mcP(\C^m)$. With the $G\times K_{G'}$ as well as $H \times K_{H'}$ multiplicity-free decomposition of $\mcH$, one allows $\mcH^{N_H \times N_{K_{G'}} }$ to be interpreted as a branching algebra from $K_{H'}$ to $K_{G'}$ as well as a branching algebra from $G$ to $H$.  This double interpretation solve two related branching problems simultaneously. Classical invariant theory
also provides a flexible approach which allows an inductive approach to the
computation of branching algebras, and makes evident natural connections
between different branching algebras.  We refer to readers to \cite{howe-tan-jeb} for more details.

\newpage

\section{Proofs}

\subsection{Proofs of the Tensor Product Formulas.}

\subsubsection{${\bf GL_n \subset GL_n \times GL_n}$}

We consider the following see-saw pair and its complexificiation:
\[
\begin{array}{ccccccc}
    U_n \times U_n & - & {\frak u}_{p,q} \oplus {\frak u}_{r,s} & \text{\tiny Complexified} & GL_n \times GL_n & - & \gl_{p,q} \oplus \gl_{r,s} \\
    \cup             &   & \cap                 & \rightsquigarrow    & \cup             &   & \cap                 \\
    U_n             & - & {\frak u}_{p+r, q+s}          &                     & GL_n             & - & \gl_{p+r,q+s}
\end{array}
\]
\vskip 10pt

Regarding the dual pair $(GL_n \times GL_n, \gl_{p,q} \oplus \gl_{r,s})$, Theorem 3.2 gives the decomposition:
\[
    \mcP\left( \left( \bbC^n \ot \bbC^p \oplus \dual{\bbC^n} \ot \bbC^q \right) \oplus
               \left( \bbC^n \ot \bbC^r \oplus \dual{\bbC^n} \ot \bbC^s \right) \right)
\]
\[  \cong \bigoplus
    \left( \F{(\mu^+, \mu^-)}{n} \ot \F{(\nu^+, \nu^-)}{n} \right) \ot
    \left( \wF{p,q}{(\mu^+, \mu^-)} \ot \wF{r,s}{(\nu^+, \nu^-)} \right)
\] where the sum is over non-negative integer partitions $\mu^+$, $\nu^+$, $\mu^-$, and
$\nu^-$ such that:
\[
\begin{array}{ll}
\ell(\mu^+) \leq p, & \ell(\mu^-) \leq q, \\
\ell(\nu^+) \leq r, & \ell(\nu^-) \leq s,  \\
\ell(\mu^+) + \ell(\mu^-) \leq n,\quad & \ell(\nu^+) + \ell(\nu^-) \leq n.
\end{array}
\]

\nt Regarding the dual pair $(GL_n, \gl_{p+r,q+s})$, Theorem 3.2 gives the decomposition:
\[
    \mcP\left( \bbC^n \ot \bbC^{p+r} \oplus \dual{\bbC^n} \ot \bbC^{q+s} \right)
    \cong \bigoplus \F{(\la^+,\la^-)}{n} \ot \wF{p+r,q+s}{(\la^+, \la^-)}
\] where the sum is over all non-negative integer partitions $\la^+$ and $\la^-$ such
that $\ell(\la^+) + \ell(\la^-) \leq n$, $\ell(\la^+) \leq p+r$ and $\ell(\la^-) \leq
q+s$.

We assume that we are in the stable range: $n \geq p+q+r+s$, so
that as a $GL_{p+r} \times GL_{q+s}$ representation (see Theorem 3.2):
\[\wF{p+r,q+s}{(\la^+, \la^-)} \cong \mcS(\bbC^{p+r} \ot \bbC^{q+s}) \otimes
\F{\la^+}{p+r} \ot \F{\la^-}{q+s}\]
As a $GL_p \times GL_q \times GL_r \times GL_s$-representation, $\wF{p+r, q+s}{(\la^+,
\la^-)}$ is equivalent to:
\[
    \mcS(\bbC^p \ot \bbC^q) \otimes \mcS(\bbC^r \ot \bbC^s) \otimes
    \mcS(\bbC^p \ot \bbC^s) \otimes \mcS(\bbC^r \ot \bbC^q) \otimes
    \F{\la^+}{p+r} \ot \F{\la^-}{q+s}
\]

\nt Note that $n \geq p+q+r+s$ implies that $n \geq p+q$ and $n \geq r+s$, so that (see Theorem 3.2):
\[\wF{p,q}{(\mu^+, \mu^-)} \cong \mcS(\bbC^p \ot \bbC^q) \otimes \F{\mu^+}{p} \ot
\F{\mu^-}{q}\]
and
\[\wF{r,s}{(\nu^+, \nu^-)} \cong \mcS(\bbC^r \ot \bbC^s) \otimes \F{\nu^+}{r} \ot
\F{\nu^-}{s}.\]

\nt Our see-saw pair implies (see (3.15)):
\[
\left[\F{(\mu^+, \mu^-)}{n} \ot \F{(\nu^+,\nu^-)}{n},
\F{(\la^+,\la^-)}{n} \right] = \left[\wF{p+r,q+s}{(\la^+,\la^-)},
\wF{p,q}{(\mu^+, \mu^-)} \ot \wF{r,s}{(\nu^+,\nu^-)}\right].
\]

\nt Using the fact that we are in the stable range:
\begin{eqnarray*}
&\left[\F{(\mu^+, \mu^-)}{n} \ot \F{(\nu^+,\nu^-)}{n},
 \F{(\la^+,\la^-)}{n} \right] \\
&=\left[\mcS(\bbC^p \ot \bbC^s) \otimes \mcS(\bbC^r \ot \bbC^q) \otimes
\F{\la^+}{p+r} \ot \F{\la^-}{q+s}, \F{\mu^+}{p} \ot \F{\mu^-}{q}
\ot \F{\nu^+}{r} \ot \F{\nu^-}{s}\right].
\end{eqnarray*}

Next we will combine the standard decompositions:
\[ \F{\la^+}{p+r} \cong \bigoplus \c{\la^+}{\al_1}{\al_2} \F{\al_1}{p} \ot \F{\al_2}{r}\]

\[ \F{\la^-}{q+s} \cong \bigoplus \c{\la^-}{\be_1}{\be_2} \F{\be_1}{q} \ot \F{\be_2}{s}\]
with the multiplicity-free decompositions (see (3.12)):
\[ \mcS( \bbC^p \ot \bbC^s) \cong \bigoplus \F{\ga_1}{p} \ot \F{\ga_1}{s} \]
\[ \mcS( \bbC^r \ot \bbC^q) \cong \bigoplus \F{\ga_2}{r} \ot \F{\ga_2}{q}. \]

\nt This implies the result:
\[
\left[\F{(\mu^+, \mu^-)}{n} \ot \F{(\nu^+,\nu^-)}{n},
\F{(\la^+,\la^-)}{n} \right] = \sum \c{\la^+}{\al_2}{\al_1}
\c{\mu^+}{\al_1}{\ga_1} \c{\nu^-}{\ga_1}{\be_2}
\c{\la^-}{\be_2}{\be_1} \c{\mu^-}{\be_1}{\ga_2}
\c{\nu^+}{\ga_2}{\al_2}.
\]

\vskip 10pt

\subsubsection{${\bf O_n \subset O_n \times O_n}$}

We consider the following see-saw pair and its complexificiation:
\vskip 5pt
\[
\begin{array}{ccccccc}
    O_n(\bbR) \times O_n(\bbR) & - & sp_{2p}(\bbR) \oplus sp_{2q}(\bbR) & \text{\tiny Complexified} & O_n \times O_n & - & \sp_{2p} \oplus \sp_{2q} \\
    \cup                       &   & \cap                               & \rightsquigarrow          & \cup                       &   & \cap                                 \\
    O_n(\bbR)                  & - & sp_{2(p+q)}(\bbR)                  &                           & O_n                  & - & \sp_{2(p+q)}
\end{array}
\]
\vskip 10pt

Regarding the dual pair $(O_n \times O_n, \sp_{2p} \oplus
\sp_{2q})$, Theorem 3.2 gives the decomposition:
\[
    \mcP\left(\bbC^n \ot \bbC^p \oplus \bbC^n \ot \bbC^q \right)
    \cong \bigoplus
    \left( \E{\mu}{n} \ot \E{\nu}{n} \right) \ot
    \left( \wE{2p}{\mu} \ot \wE{2q}{\nu} \right)
\] where the sum is over non-negative integer partitions $\mu$ and $\nu$ such that:
\[
\begin{array}{ll}
\ell(\mu) \leq p, & (\mu^\prime)_1 + (\mu^\prime)_2 \leq n, \\
\ell(\nu) \leq q, & (\nu^\prime)_1 + (\nu^\prime)_2 \leq n.
\end{array}
\]

\nt Regarding the dual pair $(O_n, \sp_{2(p+q)})$, Theorem 3.2 gives the
decomposition:
\[
    \mcP\left( \bbC^n \ot \bbC^{p+q} \right)
    \cong \bigoplus \E{\la}{n} \ot \wE{2(p+q)}{\la}
\] where the sum is over all non-negative integer partitions $\la$ such
that $\ell(\la) \leq p+q$, and $(\la^\prime)_1 +
(\la^\prime)_2 \leq n$.

We assume that we are in the stable range: $n \geq 2(p+q)$, so
that as a $GL_{p+q}$ representation (see Theorem 3.2):
\[\wE{2(p+q)}{\la} \cong \mcS(\mcS^2 \bbC^{p+q}) \otimes \F{\la}{p+q}.\]
As a $GL_p\times GL_q$-representation, $\wE{2(p+q)}{\la}$ is
equivalent to:
\[
    \mcS(\mcS^2 \bbC^p)        \otimes
    \mcS(\mcS^2 \bbC^q)        \otimes
    \mcS(\bbC^p \otimes \bbC^q) \otimes
    \F{\la}{p+q}
\]

\nt Note that $n \geq 2(p+q)$ implies that $n \geq 2p$ and $n \geq
2q$, so that (see Theorem 3.2):
\[\wE{2p}{\mu} \cong \mcS(\mcS^2 \bbC^{p}) \otimes \F{\mu}{p}\]
and
\[\wE{2q}{\nu} \cong \mcS(\mcS^2 \bbC^{q}) \otimes \F{\nu}{q}.\]

\nt Our see-saw pair implies (see (3.15)):
\[
\left[\E{\mu}{n} \ot \E{\nu}{n}, \E{\la}{n} ] = [\wE{2(p+q)}{\la},
\wE{2p}{\mu} \ot \wE{2q}{\nu}\right].
\]

\nt Using the fact that we are in the stable range: {\small
\begin{eqnarray*}
 && \left[\wE{2(p+q)}{\la}, \wE{2p}{\mu} \ot \wE{2q}{\nu}\right] \\
&&= \left[ \mcS(\mcS^2 \bbC^p)        \otimes
    \mcS(\mcS^2 \bbC^q)        \otimes
    \mcS(\bbC^p \otimes \bbC^q) \otimes
    \F{\la}{p+q}
,\mcS(\mcS^2 \bbC^{p}) \otimes \F{\mu}{p} \; \otimes \;
    \mcS(\mcS^2 \bbC^{q}) \otimes \F{\nu}{q} \right] \\
&&=\left[        \mcS(\bbC^p \otimes \bbC^q) \otimes
    \F{\la}{p+q}
,\F{\mu}{p} \otimes \F{\nu}{q} \right].
\end{eqnarray*}}
Next we will combine the decomposition:
\[ \F{\la}{p+q} \cong \bigoplus \c{\la}{\al}{\be} \F{\al}{p} \ot \F{\be}{q}\]
with the multiplicity-free decomposition (see (3.12)):
\[ \mcS( \bbC^p \ot \bbC^q) \cong \bigoplus \F{\ga}{p} \ot \F{\ga}{q} \]
to obtain the result, but first note that in the above
decompositions $\alpha$, $\beta$, and $\gamma$ range over all
non-negative integer partitions such that $\ell(\al) \leq p$,
$\ell(\be) \leq q$ and $\ell(\gamma)\leq \min(p,q)$. So we obtain:
\[
\left[\E{\mu}{n} \ot \E{\nu}{n}, \E{\la}{n}\right] = \sum_{\al,\be,\ga}
\c{\la}{\al}{\be} \c{\mu}{\al}{\ga} \c{\nu}{\be}{\ga}
\]
The above sum is over all non-negative integer partitions $\al,
\be, \ga$ such that $\ell(\al) \leq p$, $\ell(\be) \leq q$ and
$\ell(\gamma)\leq \min(p,q)$, however, the support of the
Littlewood-Richardson coefficients is contained inside the set of
such $(\al,\be,\ga)$ when we choose $p$ and $q$ such that
$\ell(\la) \leq \floor{n/2} := p+q$, with $\ell(\mu) := p$ and
$\ell(\nu) := q$.

\vskip 10pt

\subsubsection{${\bf Sp_{2n} \subset Sp_{2n} \times Sp_{2n}}$}

We consider the following see-saw pair and its complexificiation:
\vskip 5pt
\[
\begin{array}{ccccccc}
    Sp(n) \times Sp(n) & - & \so^*_{2p} \oplus \so^*_{2q} & \text{\tiny Complexified} & Sp_{2n} \times Sp_{2n} & - & \so_{2p} \oplus \so_{2q} \\
    \cup                               &   & \cap                                     & \rightsquigarrow          & \cup                   &   & \cap                     \\
    Sp(n)                      & - & so^*_{2(p+q)}                     &                           & Sp_{2n}                & - & \so_{2(p+q)}
\end{array}
\]
\vskip 10pt

Regarding the dual pair $(Sp_{2n} \times Sp_{2n}, \so_{2p} \oplus
\so_{2q})$, Theorem 3.2 gives the decomposition:
\[  \mcP\left(\bbC^{2n} \ot \bbC^p \oplus \bbC^{2n} \ot \bbC^q \right)
    \cong \bigoplus
    \left( \V{\mu}{2n} \ot \V{\nu}{2n} \right) \ot
    \left( \wV{2p}{\mu} \ot \wV{2q}{\nu} \right)
\] where the sum is over non-negative integer partitions $\mu$ and $\nu$ such that:
\[\begin{array}{ll}
\ell(\mu) \leq \min(n,p), & \ell(\nu) \leq \min(n,q).
\end{array}\]

\nt Regarding the dual pair $(Sp_{2n}, \so_{2(p+q)})$, Theorem 3.2 gives the
decomposition:
\[  \mcP\left( \bbC^{2n} \ot \bbC^{p+q} \right)
    \cong \bigoplus \V{\la}{2n} \ot \wV{2(p+q)}{\la}
\] where the sum is over all non-negative integer partitions $\la$ such
that $\ell(\la) \leq \min(n,p+q)$.

We assume that we are in the stable range: $n \geq p+q$, so
that as a $GL_{p+q}$ representation (see Theorem 3.2):
\[\wV{2(p+q)}{\la} \cong \mcS(\wedge^2 \bbC^{p+q}) \otimes \F{\la}{p+q}.\]
As a $GL_p\times GL_q$-representation, $\wV{2(p+q)}{\la}$ is
equivalent to:
\[
    \mcS(\wedge^2 \bbC^p)        \otimes
    \mcS(\wedge^2 \bbC^q)        \otimes
    \mcS(\bbC^p \otimes \bbC^q) \otimes
    \F{\la}{p+q}
\]

\nt Note that $n \geq p+q$ implies that $n \geq p$ and $n \geq q$, so that (see Theorem 3.2):
\[\wV{2p}{\mu} \cong \mcS(\wedge^2 \bbC^{p}) \otimes \F{\mu}{p}\]
and
\[\wV{2q}{\nu} \cong \mcS(\wedge^2 \bbC^{q}) \otimes \F{\nu}{q}.\]

\nt Our see-saw pair implies (see (3.15)):
\[
\left[\V{\mu}{2n} \ot \V{\nu}{2n}, \V{\la}{2n} \right] = \left[\wV{2(p+q)}{\la},
\wV{2p}{\mu} \ot \wV{2q}{\nu}\right].
\]

Using the fact that we are in the stable range: {\small
\begin{eqnarray*}
 && \left[\wV{2(p+q)}{\la}, \wV{2p}{\mu} \ot \wV{2q}{\nu}\right] \\
&&=\left[ \mcS(\wedge^2 \bbC^p)        \otimes
    \mcS(\wedge^2 \bbC^q)        \otimes
    \mcS(\bbC^p \otimes \bbC^q) \otimes
    \F{\la}{p+q}
, \mcS(\wedge^2 \bbC^{p}) \otimes \F{\mu}{p} \; \otimes \;
    \mcS(\wedge^2 \bbC^{q}) \otimes \F{\nu}{q} \right]\\
&&= \left[        \mcS(\bbC^p \otimes \bbC^q) \otimes
    \F{\la}{p+q}
,\F{\mu}{p} \otimes \F{\nu}{q}\right].
\end{eqnarray*}}
Next we will combine the decomposition:
\[ \F{\la}{p+q} \cong \bigoplus \c{\la}{\al}{\be} \F{\al}{p} \ot \F{\be}{q}\]
with the multiplicity-free decomposition (see (3.12)):
\[ \mcS( \bbC^p \ot \bbC^q) \cong \bigoplus \F{\ga}{p} \ot \F{\ga}{q} \]
to obtain the result, but first note that in the above
decompositions $\alpha$, $\beta$, and $\gamma$ range over all
non-negative integer partitions such that $\ell(\al) \leq p$,
$\ell(\be) \leq q$ and $\ell(\gamma)\leq \min(p,q)$. So we obtain:
\[
\left[ \V{\mu}{2n} \ot \V{\nu}{2n}, \V{\la}{2n} \right] = \sum_{\al,\be,\ga}
\c{\la}{\al}{\be} \c{\mu}{\al}{\ga} \c{\nu}{\be}{\ga}
\]
The above sum is over all non-negative integer partitions $\al,
\be, \ga$ such that $\ell(\al) \leq p$, $\ell(\be) \leq q$ and
$\ell(\gamma)\leq \min(p,q)$, however, the support of the
Littlewood-Richardson coefficients is contained inside the set of
such $(\al,\be,\ga)$ when we choose $p$ and $q$ such that
$\ell(\la) \leq \floor{n/2} := p+q$, with $\ell(\mu) := p$ and $\ell(\nu) := q$.

\newpage

\subsection{Proofs of the Direct Sum Branching Rules.}

\subsubsection{$ {\bf GL_n    \times GL_m    \subset GL_{n+m} }   $}

We consider the following see-saw pair and its complexificiation:
\vskip 5pt
\[
\begin{array}{ccccccc}
    U_{n+m}            & - & u_{p,q}                & \text{\tiny Complexified} & GL_{n+m}         & - & \gl_{p,q} \\
    \cup              &   & \cap                   & \rightsquigarrow          & \cup             &   & \cap                       \\
    U_n \times U_m  & - & u_{p,q} \oplus u_{p,q} &                           & GL_n \times GL_m & - & \gl_{p,q} \oplus \gl_{p,q}
\end{array}
\]
\vskip 10pt

Regarding the dual pair $(GL_{n+m}, \gl_{p+q})$, Theorem 3.2 gives the decomposition:
\[  \mcP\left(\bbC^{n+m} \ot \bbC^p \oplus \dual{\bbC^{n+m}} \ot \bbC^q \right)
    \cong \bigoplus \F{(\la^+, \la^-)}{n+m} \ot \wF{p,q}{(\la^+,\la^-)}
\] where the sum is over non-negative integer partitions $\la^+$ and $\la^-$ such that $\ell(\la^+) \leq p$,
$\ell(\la^-) \leq q$ and $\ell(\la^+) + \ell(\la^-) \leq n+m$.
Regarding the dual pair $(GL_n \times GL_m , \gl_{p+q} \oplus \gl_{p+q})$, Theorem 3.2 gives the
decomposition:
\[
\begin{array}{cc}
& \mcP\left(
      \bbC^n  \ot \bbC^p \oplus
\dual{\bbC^n} \ot \bbC^q \oplus
      \bbC^m  \ot \bbC^p \oplus
\dual{\bbC^m} \ot \bbC^q
\right)\\
& \cong \bigoplus
\left(
    \F{(\mu^+, \mu^-)}{n} \ot \F{(\nu^+, \nu^-)}{m}
\right) \ot
\left(
    \wF{p,q}{(\mu^+, \mu^-)} \ot \wF{p,q}{(\nu^+, \nu^-)}
\right)
\end{array}
\]
where the sum is over all non-negative integer partitions $\mu^+$,
$\mu^-$, $\nu^+$ and $\nu^-$ such that:
\[
\begin{array}{ll}
\ell(\mu^+) + \ell(\mu^-) \leq n, & \ell(\nu^+) + \ell(\nu^-) \leq m, \\
\ell(\mu^+) \leq p,\quad              & \ell(\mu^-) \leq q, \\
\ell(\nu^+) \leq p,               & \ell(\nu^-) \leq q.
\end{array}
\]
We assume that we are in the stable range: $\min(n,m) \geq p+q$, so
that as a $GL_p \times GL_q$ representation (see Theorem 3.2):
\[ \wF{p,q}{(\mu^+, \mu^-)} \cong \mcS(\bbC^p \ot \bbC^q) \ot \F{\mu^+}{p} \ot \F{\mu^-}{q} \]
\[ \wF{p,q}{(\nu^+, \nu^-)} \cong \mcS(\bbC^p \ot \bbC^q) \ot \F{\nu^+}{p} \ot \F{\nu^-}{q} \]

\nt Note that $\min(n,m) \geq p+q$ implies that $n+m \geq p+q$, so that (see Theorem 3.2):
\[ \wF{p,q}{(\la^+, \la^-)} \cong \mcS(\bbC^p \ot \bbC^q) \ot \F{\la^+}{p} \ot \F{\la^-}{q} \]
Our see-saw pair implies (see (3.15)):
\[
\left[\wF{p,q}{(\mu^+, \mu^-)}\ot\wF{p,q}{(\nu^+, \nu^-)},\wF{p,q}{(\la^+, \la^-)}\right] =
\left[\F{(\la^+, \la^-)}{n+m},\F{(\mu^+,\mu^-)}{n} \ot \F{(\nu^+,\nu^-)}{m}\right].
\]
Using the fact that we are in the stable range:{\small
\begin{eqnarray*}
&& \left[\wF{p,q}{(\mu^+, \mu^-)}\ot\wF{p,q}{(\nu^+, \nu^-)},\wF{p,q}{(\la^+, \la^-)}\right] \\
&&= \left[ \left( \mcS(\bbC^p \ot \bbC^q) \otimes \F{\mu^+}{p} \ot \F{\mu^-}{q} \right) \otimes
    \left( \mcS(\bbC^p \ot \bbC^q) \otimes \F{\nu^+}{p} \ot \F{\nu^-}{q} \right)
,\mcS(\bbC^p \ot \bbC^q) \otimes \F{\la^+}{p} \ot \F{\la^-}{q} \right] \\
&&= \left[\mcS(\bbC^p \ot \bbC^q) \otimes \F{\mu^+}{p} \ot \F{\mu^-}{q} \otimes
                                    \F{\nu^+}{p} \ot \F{\nu^-}{q}
,\F{\la^+}{p} \ot \F{\la^-}{q}\right].
\end{eqnarray*}}
Next combine the above decomposition with (see (3.12)):
\[ \mcS( \bbC^p \ot \bbC^q) \cong \bigoplus \F{\de}{p} \ot \F{\de}{q} \]
where the sum is over all non-negative integer partitions $\de$ with at most $\min(p,q)$ parts.
We then obtain:
\[
\begin{array}{cc}
& \left[\F{(\la^+, \la^-)}{n+m},\F{(\mu^+,\mu^-)}{n} \ot \F{(\nu^+,\nu^-)}{m}\right] \\
&= \left[ \left( \bigoplus_\de \F{\de  }{p} \ot \F{\de  }{q} \right) \otimes
  \left(           \F{\mu^+}{p} \ot \F{\mu^-}{q} \right) \otimes
  \left(           \F{\nu^+}{p} \ot \F{\nu^-}{q} \right),\F{\la^+}{p} \ot \F{\la^-}{q} \right].
\end{array}
\]
We combine this fact with the following two tensor product
decompositions:
\[
\begin{array}{lllllll}
 \F{\mu^+}{p} \otimes \F{\nu^+}{p} & \cong & \bigoplus \c{\ga^+}{\mu^+}{\nu^+} \F{\ga^+}{p}
& \text{and } & \F{\mu^-}{q} \otimes \F{\nu^-}{q} & \cong & \bigoplus \c{\ga^-}{\mu^-}{\nu^-} \F{\ga^-}{q}
\end{array}
\]
(where $\ga^+$ and $\ga^-$ have at most $p$ and $q$ parts respectively) and then tensor the constituents
with $\F{\de}{p} \ot \F{\de}{q}$,
\[
\begin{array}{lllllll}
 \F{\ga^+}{p} \otimes \F{\de}{p}  & \cong & \bigoplus \c{\la^+}{\ga^+}{\de} \F{\la^+}{p}
& \text{and} & \F{\ga^-}{q} \otimes \F{\de}{q} & \cong & \bigoplus \c{\la^-}{\ga^-}{\de} \F{\la^-}{q}
\end{array}
\] to obtain the result.

\vskip 10pt

\subsubsection{$ {\bf O_n    \times O_m     \subset  O_{n+m} }  $}

We consider the following see-saw pair and its complexificiation:
\vskip 5pt
\[
\begin{array}{ccccccc}
    O_{n+m}(\R)            & - & \sp_{2k}(\R)                 & \text{\tiny Complexified} & O_{n+m}          & - & \sp_{2k} \\
    \cup              &   & \cap                     & \rightsquigarrow          & \cup             &   & \cap                       \\
    O_n(\R) \times O_m(\R)  & - & \sp_{2k}(\R) \oplus \sp_{2k}(\R) &                           & O_n \times O_m   & - & \sp_{2k} \oplus \sp_{2k}
\end{array}
\]
\vskip 10pt

Regarding the dual pair $(O_{n+m}, \sp_{2k})$, Theorem 3.2 gives the decomposition:
\[  \mcP\left(\bbC^{n+m} \ot \bbC^k \right)
    \cong \bigoplus \E{\la}{n+m} \ot \wE{2k}{\la}
\] where the sum is over non-negative integer partitions $\la$ such that $\ell(\la) \leq k$ and $(\la^\prime)_1 + (\la^\prime)_2 \leq n+m$.  Regarding the dual pair $(O_n \times O_m, \sp_{2k} \oplus \sp_{2k})$, Theorem 3.2 gives the
decomposition:
\[
\mcP\left(
      \bbC^n  \ot \bbC^k \oplus
      \bbC^m  \ot \bbC^k
\right) \cong \bigoplus
\left(  \E{\mu}{n}  \ot \E{\nu}{m}    \right) \ot
\left( \wE{2k}{\mu} \ot \wE{2k}{\nu} \right)
\] where the sum is over all non-negative integer partitions $\mu$ and $\nu$ such that
$\ell(\mu) \leq k$, $\ell(\nu) \leq k$, $(\mu^\prime)_1 + (\mu^\prime)_2 \leq n$ and $(\nu^\prime)_1 + (\nu^\prime)_2 \leq m$.

We assume that we are in the stable range: $\min(n,m) \geq 2k$, so
that as $GL_k$ representations (see Theorem 3.2):
\[
\begin{array}{lll}
\wE{2k}{\mu} \cong \mcS(\mcS^2 \bbC^k ) \ot \F{\mu}{k}
& \text{ and } &
\wE{2k}{\nu} \cong \mcS(\mcS^2 \bbC^k ) \ot \F{\nu}{k}.
\end{array}
\]

\nt Note that $\min(n,m) \geq 2k$ implies that $n + m \geq 2k$, so that (see Theorem 3.2):
\[ \wE{2k}{\la} \cong \mcS(\mcS^2 \bbC^k ) \ot \F{\la}{k} \]

Our see-saw pair implies (see (3.15)):
\[
\left[\wE{2k}{\mu}\ot\wE{2k}{\nu},\wE{2k}{\la}\right] =
\left[\E{\la}{n+m},\E{\mu}{n} \ot \F{\nu}{m}\right].
\]
Using the fact that we are in the stable range:{\small
\begin{eqnarray*}
\left[\wE{2k}{\mu}\ot\wE{2k}{\nu},\wE{2k}{\la}\right]
&&= \left[ \left( \mcS(\mcS^2\bbC^k) \otimes \F{\mu}{k}  \right) \otimes
      \left( \mcS(\mcS^2\bbC^k) \otimes \F{\nu}{k}  \right)
    ,\mcS(\mcS^2\bbC^k) \otimes \F{\la}{k} \right] \\
&&= \left[ \mcS(\mcS^2\bbC^k) \otimes \F{\mu}{k} \otimes \F{\nu}{k},\F{\la}{k} \right].
\end{eqnarray*}}
Next combine with the well-known multiplicity-free decomposition (see for instance, Theorem 3.1 of \cite{howe-schur} on page 32):
\[ \mcS( \mcS^2 \bbC^k ) \cong \bigoplus \F{2\de}{k} \]
where the sum is over all non-negative integer partitions $\de$ with at most $k$ parts.
We then obtain:
\[
\left[\E{\la}{n+m},\E{\mu}{n} \ot \E{\nu}{m}\right] =
\left[ \left(
\bigoplus_\de \F{2\de  }{k} \right) \otimes
       \F{\mu}{k} \otimes \F{\nu}{k}, \F{\la}{k} \right].
\]
Combine this fact with the following two tensor product decompositions:
\[
\begin{array}{lllllll}
 \F{\mu}{k} \otimes \F{\nu}{k} & \cong & \bigoplus \c{\ga}{\mu}{\nu} \F{\ga}{k}
& \text{and } & \F{\ga}{k} \otimes \F{2\de}{k} & \cong & \bigoplus \c{\la}{\ga}{2 \de} \F{\la}{k}
\end{array}
\]
(where $\ga$ is a non-negative integer partition with at most $k$ parts) and the result follows.

\subsubsection{$ {\bf Sp_{2n} \times Sp_{2m} \subset Sp_{2(n+m)}} $}

We consider the following see-saw pair and its complexificiation:
\vskip 10pt
\[
\begin{array}{ccccccc}
    Sp(n+m)             & - & \so^*_{2k}                   & \text{\tiny Complexified} & Sp_{2(n+m)}            & - & \so_{2k} \\
    \cup                &   & \cap                         & \rightsquigarrow          & \cup                   &   & \cap                       \\
    Sp(n) \times Sp(m)  & - & \so^*_{2k} \oplus \so^*_{2k} &                           & Sp_{2n} \times Sp_{2m} & - & \so_{2k} \oplus \so_{2k}
\end{array}
\]
\vskip 10pt

Regarding the dual pair $(Sp_{2(n+m)},\so_{2k})$, Theorem 3.2 gives the decomposition:
\[  \mcP\left(\bbC^{2(n+m)} \ot \bbC^k \right)
    \cong \bigoplus \V{\la}{2(n+m)} \ot \wV{2k}{\la}
\] where the sum is over non-negative integer partitions $\la$ such that\linebreak  $\ell(\la) \leq \min(n+m,k)$. Regarding the dual pair $(Sp_{2n} \times Sp_{2m}, \so_{2k} \oplus \so_{2k})$, Theorem 3.2 gives the
decomposition:
\[
\mcP\left(
      \bbC^{2n}  \ot \bbC^k \oplus
      \bbC^{2m}  \ot \bbC^k
\right) \cong \bigoplus
\left(  \V{\mu}{2n}  \ot \V{\nu}{2m}    \right) \ot
\left( \wV{2k}{\mu} \ot \wV{2k}{\nu} \right)
\] where the sum is over all non-negative integer partitions $\mu$ and $\nu$ such that
$\ell(\mu) \leq \min(n,k)$, $\ell(\nu) \leq \min(n,k)$.

We assume that we are in the stable range: $\min(n,m) \geq k$, so that as $GL_k$ representations (see Theorem 3.2):
\[
\begin{array}{lll}
\wV{2k}{\mu} \cong \mcS(\wedge^2 \bbC^k ) \ot \F{\mu}{k}
& \text{ and } &
\wV{2k}{\nu} \cong \mcS(\wedge^2 \bbC^k ) \ot \F{\nu}{k}.
\end{array}
\]

\nt Note that $\min(n,m) \geq k$ implies that $n + m \geq k$, so that (see Theorem 3.2):
\[ \wV{2k}{\la} \cong \mcS(\wedge^2 \bbC^k ) \ot \F{\la}{k} \]
Our see-saw pair implies (see (3.15)):
\[
\left[\wV{2k}{\mu}\ot\wV{2k}{\nu},\wV{2k}{\la}\right] =
\left[\V{\la}{n+m},\V{\mu}{n} \ot \V{\nu}{m}\right].
\]
Using the fact that we are in the stable range: {\small
\begin{eqnarray*}
\left[\wV{2k}{\mu}\ot\wV{2k}{\nu},\wV{2k}{\la}\right]
&=& \left[ \left( \mcS(\wedge^2\bbC^k) \otimes \F{\mu}{k}  \right) \otimes
      \left( \mcS(\wedge^2\bbC^k) \otimes \F{\nu}{k}  \right)
    ,\mcS(\wedge^2\bbC^k) \otimes \F{\la}{k} \right] \\
&=&\left[ \mcS(\wedge^2\bbC^k) \otimes \F{\mu}{k} \otimes \F{\nu}{k},\F{\la}{k} \right].
\end{eqnarray*}}
Next combine with the well-known multiplicity-free decomposition (see for instance, Theorem 3.8.1 of \cite{howe-schur} on page 44):
\[ \mcS( \wedge^2 \bbC^k ) \cong \bigoplus \F{(2\de)^\prime}{k} \]
where the sum is over all non-negative integer partitions $\de$ such that $(2\de)^\prime$ has at most $k$ parts.
We then obtain:
\[
\left[\V{\la}{2(n+m)},\V{\mu}{2n} \ot \V{\nu}{2m}\right] =
\left[ \left( \bigoplus_\de \F{(2\de)^\prime  }{k} \right) \otimes
       \F{\mu}{k} \otimes \F{\nu}{k}, \F{\la}{k} \right].
\]
Combine this fact with the following two tensor product decompositions:
\[
\begin{array}{lllllll}
 \F{\mu}{k} \otimes \F{\nu}{k} & \cong & \bigoplus \c{\ga}{\mu}{\nu} \F{\ga}{k}
& \text{and } & \F{\ga}{k} \otimes \F{(2\de)^\prime}{k} & \cong & \bigoplus \c{\la}{\ga}{(2 \de)^\prime} \F{\la}{k}
\end{array}
\]
(where $\ga$ is a non-negative integer partition with at most $k$ parts) and the result follows.

\vskip 10pt

\subsection{Proofs of the Polarization Branching Rules.}

\subsubsection{$ {\bf GL_n \subset  O_{2n} }$}

We consider the following see-saw pair and its complexificiation:
\vskip 5pt
\[
\begin{array}{ccccccc}
    O_{2n}(\bbR) & - & \sp_{2k}(\bbR) & \text{\tiny Complexified} & O_{2n} & - & \sp_{2k} \\
    \cup  &   & \cap           & \rightsquigarrow          & \cup   &   & \cap     \\
    U(n)  & - & u_{k,k}      &                           & GL_n   & - & \gl_{k,k}
\end{array}
\]
\vskip 10pt

Regarding the dual pair $(O_{2n}, \sp_{2k} )$, Theorem 3.2 gives the
decomposition:
\[
\mcP\left( \bbC^{2n} \ot \bbC^k \right) \cong \bigoplus \E{\la}{2n} \ot \wE{2k}{\la}
\] where the sum is over all non-negative integer partitions $\la$ such that
$\ell(\la) \leq k$ and $(\la^\prime)_1 + (\la^\prime)_2 \leq 2n$. Since the standard $O_{2n}$ representation $\C^{2n} \simeq \C^n \oplus (\C^n)^\ast$ as a $GL_n$ representation,
regarding the dual pair $(GL_{n} , \gl_{k,k})$, Theorem 3.2 gives the decomposition:
\[
\mcP \left(
            \bbC^n  \ot \bbC^k  \otimes
      \dual{\bbC^n} \ot \bbC^k \right)
    \cong \bigoplus \F{(\mu^+,\mu^-)}{n} \ot \wF{k,k}{(\mu^+, \mu^-)}
\] where the sum is over non-negative integer partitions $\mu^+$ and $\mu^-$ with at most $k$
parts such that $\ell(\mu^+) + \ell(\mu^-) \leq n$.

We assume that we are in the stable range: $n \geq 2k$, so that as a $GL_k \times GL_k$ representation (see Theorem 3.2):
\[ \wF{k,k}{(\mu^+,\mu^-)} \cong \mcS(\bbC^k \otimes \bbC^k) \otimes \F{\mu^+}{k} \otimes \F{\mu^-}{k}. \]

\nt Note that $n \geq 2k$ implies that $n \geq k$, so that as $GL_k$ representations (see Theorem 3.2):
\[ \wE{2k}{\la} \cong \mcS(\mcS^2 \bbC^k ) \ot \F{\la}{k} \]
Our see-saw pair implies (see (3.15)):
\[
\left[\wF{k,k}{(\mu^+,\mu^-)},\wE{2k}{\la} \right] =
\left[\E{\la}{2n},\F{(\mu^+,\mu^-)}{n}\right].
\]
Using the fact that we are in the stable range:
\begin{eqnarray*}
\left[\wF{k,k}{(\mu^+,\mu^-)},\wE{2k}{\la} \right]
&=& \left[     \mcS(\bbC^k \otimes \bbC^k) \otimes \F{\mu^+}{k} \otimes \F{\mu^-}{k},
               \mcS(\mcS^2 \bbC^k ) \ot \F{\la}{k}
\right] \\
&=& \left[ \mcS(\wedge^2 \bbC^k) \otimes \F{\mu^+}{k} \otimes \F{\mu^-}{k},\F{\la}{k} \right].
\end{eqnarray*}
Note that in the above we used the fact that as a $GL_k$-representation,\linebreak
$\otimes^2 \bbC^k \cong \mcS^2 \bbC^k \oplus \wedge^2 \bbC^k$.

Next combine with the well-known multiplicity-free decomposition:
\[ \mcS( \wedge^2 \bbC^k ) \cong \bigoplus \F{(2\de)^\prime}{k} \]
where the sum is over all non-negative integer partitions $\de$ such that $(2\de)^\prime$ has at most $k$ parts.
We then obtain:
\[
\left[\E{\la}{2n},\F{(\mu^+,\mu^-)}{n}\right] =
\left[ \left( \bigoplus_\de \F{(2\de)^\prime  }{k} \right) \otimes
       \F{\mu^+}{k} \otimes \F{\mu^-}{k}, \F{\la}{k} \right].
\]
Combine this fact with the following two tensor product decompositions:
\[
\begin{array}{lllllll}
 \F{\mu^+}{k} \otimes \F{\mu^-}{k} & \cong & \bigoplus \c{\ga}{\mu^+}{\mu^-} \F{\ga}{k}
& \text{and } & \F{\ga}{k} \otimes \F{(2\de)^\prime}{k} & \cong & \bigoplus \c{\la}{\ga}{(2 \de)^\prime} \F{\la}{k}
\end{array}
\]
(where $\ga$ is a non-negative integer partition with at most $k$ parts) and the result follows.

\vskip 10pt

\subsubsection{$ {\bf GL_n \subset Sp_{2n} }$}

We consider the following see-saw pair and its complexificiation:
\vskip 5pt
\[
\begin{array}{ccccccc}
    Sp(n) & - & \so^*_{2k} & \text{\tiny Complexified} & Sp_{2n} & - & \so_{2k} \\
    \cup  &   & \cap     & \rightsquigarrow          & \cup   &   & \cap     \\
    U(n)  & - & u_{k,k}   &                           & GL_n   & - & \gl_{k,k}
\end{array}
\]
\vskip 10pt

Regarding the dual pair $(Sp_{2n}, \so_{2k} )$, Theorem 3.2 gives the
decomposition:
\[
\mcP\left( \bbC^{2n} \ot \bbC^k \right) \cong \bigoplus \V{\la}{2n} \ot \wV{2k}{\la}
\] where the sum is over all non-negative integer partitions $\la$ such that
$\ell(\la) \leq \min(n,k)$. Since the standard $Sp_{2n}$ representation $\C^{2n} \simeq \C^n \oplus (\C^n)^\ast$ as a $GL_n$ representation, regarding the dual pair $(GL_{n} ,\gl_{k,k})$, Theorem 3.2 gives the decomposition:
\[
\mcP \left(
             \bbC^n  \ot \bbC^k \otimes
       \dual{\bbC^n} \ot \bbC^k \right)
    \cong \bigoplus \F{(\mu^+,\mu^-)}{n} \ot \wF{k,k}{(\mu^+, \mu^-)}
\] where the sum is over non-negative integer partitions $\mu^+$ and $\mu^-$ with at most $k$
parts such that $\ell(\mu^+) + \ell(\mu^-) \leq n$.

We assume that we are in the stable range: $n \geq 2k$, so that as $GL_k \times GL_k$ representations (see Theorem 3.2):
\[ \wF{k,k}{(\mu^+,\mu^-)} \cong \mcS(\bbC^k \otimes \bbC^k) \otimes \F{\mu^+}{k} \otimes \F{\mu^-}{k}. \]

\nt Note that $n \geq 2k$ implies that $n \geq k$, so that as a $GL_k$ representation (see Theorem 3.2):
\[ \wV{2k}{\la} \cong \mcS(\wedge^2 \bbC^k ) \ot \F{\la}{k} \]
Our see-saw pair implies (see (3.15)):
\[
\left[\wF{k,k}{(\mu^+,\mu^-)},\wV{2k}{\la} \right] =
\left[\V{\la}{2n},\F{(\mu^+,\mu^-)}{n}\right].
\]
Using the fact that we are in the stable range:
\begin{eqnarray*}
\left[\wF{k,k}{(\mu^+,\mu^-)},\wE{2k}{\la} \right]
&=& \left[     \mcS(\bbC^k \otimes \bbC^k) \otimes \F{\mu^+}{k} \otimes \F{\mu^-}{k},
               \mcS(\wedge^2 \bbC^k ) \ot \F{\la}{k}
\right] \\
&=& \left[ \mcS(\mcS^2 \bbC^k) \otimes \F{\mu^+}{k} \otimes \F{\mu^-}{k},\F{\la}{k} \right].
\end{eqnarray*}

Note that in the above we used the fact that as a $GL_k$-representation,\linebreak
$\otimes^2 \bbC^k \cong \mcS^2 \bbC^k \oplus \wedge^2 \bbC^k$.

Next combine with the well-known multiplicity-free decomposition:
\[ \mcS( \mcS^2 \bbC^k ) \cong \bigoplus \F{2\de}{k} \]
where the sum is over all non-negative integer partitions $\de$ with at most $k$ parts.
We then obtain:
\[
\left[\V{\la}{2n},\F{(\mu^+,\mu^-)}{n}\right] =
\left[ \left( \bigoplus_\de \F{2\de}{k} \right) \otimes
       \F{\mu^+}{k} \otimes \F{\mu^-}{k}, \F{\la}{k} \right].
\]
Combine this fact with the following two tensor product decompositions:
\[
\begin{array}{lllllll}
 \F{\mu^+}{k} \otimes \F{\mu^-}{k} & \cong & \bigoplus \c{\ga}{\mu^+}{\mu^-} \F{\ga}{k}
& \text{and } & \F{\ga}{k} \otimes \F{2\de}{k} & \cong & \bigoplus \c{\la}{\ga}{2 \de} \F{\la}{k}
\end{array}
\]
(where $\ga$ is a non-negative integer partition with at most $k$ parts) and the result follows.

\subsection{Proofs of the Bilinear Form Branching Rules}

\subsubsection{$  {\bf O_n    \subset GL_n }   $}

We consider the following see-saw pair and its complexificiation:
\vskip 5pt
\[
\begin{array}{ccccccc}
    GL_n(\R)  & - &GL_{p+q}(\R)              & \text{\tiny Complexified} & GL_n & - & \gl_{p,q} \\
    \cup  &   & \cap                & \rightsquigarrow          & \cup &   & \cap   \\
    O_n(\R)  & - & \sp_{2(p+q)}(\bbR) &                           & O_n  & - & \sp_{2(p+q)}
\end{array}
\]
\vskip 10pt

Regarding the dual pair $(GL_n, \gl_{p,q})$, Theorem 3.2 gives the decomposition:
\[
\mcP \left(
            \bbC^n  \ot \bbC^p \otimes
      \dual{\bbC^n} \ot \bbC^q \right)
    \cong \bigoplus \F{(\la^+,\la^-)}{n} \ot \wF{p,q}{(\la^+, \la^-)}
\] where the sum is over non-negative integer partitions $\la^+$ and $\la^-$ with at most $p$ and $q$ parts respectively and such that $\ell(\la^+) + \ell(\la^-) \leq n$.  Noting that $\C^n \simeq (\C^n)^\ast$ as an $O_n$ representation, regarding the dual pair $(O_n, \sp_{2(p+q)} )$, Theorem 3.2 gives the decomposition:
\[
\mcP\left( \bbC^n \ot \bbC^{p+q} \right) \cong \bigoplus \E{\la}{n} \ot \wE{2(p+q)}{\la}
\] where the sum is over all non-negative integer partitions $\la$ such that
$\ell(\la) \leq p+q$ and $(\la^\prime)_1 + (\la^\prime)_2 \leq n$.

We assume that we are in the stable range: $n \geq 2(p+q)$, so that as $GL_{p+q}$ representations (see Theorem 3.2):
\[ \wE{2(p+q)}{\mu} \cong \mcS(\mcS^2 \bbC^{p+q} ) \ot \F{\mu}{p+q}. \]

\nt Note that $n \geq 2(p+q)$ implies that $n \geq p+q$, so that as $GL_p \times GL_q$ representations (see Theorem 3.2):
\[ \wF{p,q}{(\la^+,\la^-)} \cong \mcS(\bbC^p \ot \bbC^q) \otimes \F{\la^+}{p} \ot \F{\la^-}{q}. \]

Our see-saw pair implies (see (3.15):
\[
\left[\wE{2(p+q)}{\mu},\wF{p,q}{(\la^+,\la^-)}\right] =
\left[\F{(\la^+,\la^-)}{n},\E{\mu}{n}\right].
\]
Using the fact that we are in the stable range:
\begin{eqnarray*}
&& \left[\wE{2(p+q)}{\mu},\wF{p,q}{(\la^+,\la^-)}\right] \\
&&=
\left[
\mcS(\mcS^2 \bbC^{p+q} ) \ot \F{\mu}{p+q},
\mcS(\bbC^p \otimes \bbC^q) \otimes \F{\la^+}{p} \otimes \F{\la^-}{q} \right] \\
&&=
\left[ \mcS(\mcS^2 \bbC^p \oplus \mcS^2 \bbC^q \oplus \bbC^p \ot \bbC^q ) \ot \F{\mu}{p+q},
\mcS(\bbC^p \ot \bbC^q) \otimes \F{\la^+}{p} \ot \F{\la^-}{q}
\right]\\
&&=
\left[\mcS(\mcS^2 \bbC^p)  \ot \mcS( \mcS^2 \bbC^q) \otimes \F{\mu}{p+q},
\F{\la^+}{p} \ot \F{\la^-}{q} \right]
\end{eqnarray*}

Next we combine the decompositions:
\[\F{\mu}{p+q} \cong \bigoplus \c{\mu}{\al}{\be} \F{\al}{p}\ot \F{\be}{q}, \]
with the multiplicity-free decompositions:
\[\mcS(\mcS^2 \bbC^p ) \cong \bigoplus \F{2\ga}{p}\qquad  \text{and} \qquad
\mcS(\mcS^2 \bbC^q ) \cong \bigoplus \F{2\de}{q} \]
where the sums are over all non-negative integer partitions $\ga$ and $\de$ with
at most $p$ and $q$ parts respectively. We then obtain:
\[
\left[\F{(\la^+,\la^-)}{n},\E{\mu}{n}\right] =
\left[
\left( \bigoplus \F{2\ga}{p} \right) \ot
\left( \bigoplus \F{2\de}{q} \right) \otimes
\left( \bigoplus \c{\mu}{\al}{\be} \F{\al}{p}\ot \F{\be}{q} \right),
\F{\la^+}{p} \ot \F{\la^-}{q}\right].
\]
Combine this fact with the following two tensor product decompositions:
\[
\begin{array}{lllllll}
 \F{\al}{p} \otimes \F{2\ga}{p} & \cong & \bigoplus \c{\la^+}{\al}{2\ga} \F{\la^+}{p}
& \text{and } &
 \F{\be}{q} \otimes \F{2\de}{q} & \cong & \bigoplus \c{\la^-}{\be}{2\de} \F{\la^-}{q}
\end{array}
\]
and the result follows.

\vskip 10pt

\subsubsection{$ {\bf Sp_{2n} \subset GL_{2n} }$}

We consider the following see-saw pair and its complexificiation:
\vskip 5pt
\[
\begin{array}{ccccccc}
    U_{2n}  & - & {\frak u}_{p,q}         & \text{\tiny Complexified} & GL_{2n} & - & \gl_{p,q} \\
     \cup  &   & \cap           & \rightsquigarrow          & \cup    &   & \cap   \\
    Sp(n)  & - & \so^*_{2(p+q)} &                           & Sp_{2n} & - & \so_{2(p+q)}
\end{array}
\]
\vskip 10pt

Regarding the dual pair $(GL_{2n}, \gl_{p,q})$, Theorem 3.2 gives the decomposition:
\[
\mcP \left(
            \bbC^{2n}  \ot \bbC^p  \otimes
     \dual{\bbC^{2n}}  \ot \bbC^q  \right)
    \cong \bigoplus \F{(\la^+,\la^-)}{2n} \ot \wF{p,q}{(\la^+, \la^-)}
\] where the sum is over non-negative integer partitions $\la^+$ and $\la^-$ with at most $p$ and $q$ parts respectively and such that $\ell(\la^+) + \ell(\la^-) \leq 2n$.  Noting that $(\C^{2n})^\ast \simeq \C^{2n}$ as $Sp_{2n}$ modules, regarding the dual pair $(Sp_{2n}, \so_{2(p+q)} )$, Theorem 3.2 gives the decomposition:
\[
\mcP\left( \bbC^{2n} \ot \bbC^{p+q} \right) \cong \bigoplus \V{\la}{2n} \ot \wV{2(p+q)}{\la}
\] where the sum is over all non-negative integer partitions $\la$ such that
$\ell(\la) \leq \min(2n,p+q)$.

We assume that we are in the stable range: $n \geq p+q$, so that as $GL_{p+q}$ representations (see Theorem 3.2):
\[ \wV{2(p+q)}{\mu} \cong \mcS(\wedge^2 \bbC^{p+q} ) \otimes \F{\mu}{p+q} \]

\nt Note that $n \geq p+q$ implies that $2n \geq p+q$, so that as $GL_p \times GL_q$ representations (see Theorem 3.2):
\[ \wF{p,q}{(\la^+,\la^-)} \cong \mcS(\bbC^p \ot \bbC^q) \otimes \F{\la^+}{p} \ot \F{\la^-}{q} \]

Our see-saw pair implies (see (3.15)):
\[
\left[\wV{2(p+q)}{\mu},\wF{p,q}{(\la^+,\la^-)}\right] =
\left[\F{(\la^+,\la^-)}{2n},\V{\mu}{2n}\right].
\]
Using the fact that we are in the stable range:
\begin{eqnarray*}
&& \left[\wV{2(p+q)}{\mu},\wF{p,q}{(\la^+,\la^-)}\right] \\
&&= \left[ \mcS(\wedge^2 \bbC^{p+q} ) \otimes
\F{\mu}{p+q},
\mcS(\bbC^p \ot \bbC^q) \otimes \F{\la^+}{p} \ot \F{\la^-}{q} \right] \\
&&= \left[ \mcS(\wedge^2 \bbC^p \oplus \wedge^2 \bbC^q \oplus \bbC^p \ot \bbC^q ) \otimes
\F{\mu}{p+q}, \mcS(\bbC^p \ot \bbC^q) \otimes \F{\la^+}{p} \ot \F{\la^-}{q}
\right]\\
&&= \left[\mcS(\wedge^2 \bbC^p)  \ot \mcS( \wedge^2 \bbC^q) \otimes \F{\mu}{p+q},
\F{\la^+}{p} \ot \F{\la^-}{q} \right]
\end{eqnarray*}

Next we combine the decompositions:
\[\F{\mu}{p+q} \cong \bigoplus \c{\mu}{\al}{\be} \F{\al}{p}\ot \F{\be}{q}, \]
with the multiplicity-free decompositions:
\[\mcS(\wedge^2 \bbC^p ) \cong \bigoplus \F{(2\ga)^\prime}{p} \qquad \text{and}\qquad \mcS(\wedge^2 \bbC^q ) \cong \bigoplus \F{(2\de)^\prime}{q}, \]
where the sums are over all non-negative integer partitions $\ga$
and $\de$ such that $(2\ga)^\prime$ and $(2\de)^\prime$ have at
most $p$ and $q$ parts respectively. We then obtain: {\small
\[
\left[\F{(\la^+,\la^-)}{2n},\V{\mu}{2n}\right] = \left[ \left( \bigoplus \F{(2\ga)^\prime}{p} \right) \ot
\left( \bigoplus \F{(2\de)^\prime}{q} \right) \otimes \left( \bigoplus \c{\mu}{\al}{\be}
\F{\al}{p}\ot \F{\be}{q} \right), \F{\la^+}{p} \ot \F{\la^-}{q} \right].
\]}
Combine with the following two tensor product decompositions:
\[
\begin{array}{lllllll}
 \F{\al}{p} \otimes \F{(2\ga)^\prime}{p} & \cong & \bigoplus \c{\la^+}{\al}{(2\ga)^\prime} \F{\la^+}{p}
& \text{and } &
 \F{\be}{q} \otimes \F{(2\de)^\prime}{q} & \cong & \bigoplus \c{\la^-}{\be}{(2\de)^\prime} \F{\la^-}{q}
\end{array}
\]
and the result follows.

\end{document}